%
% Flatness of noetherian Hopf algebras over coideal subalgebras
%	S. Skryabin
%
%	Plain TeX + AMS fonts ( amssym.def )
%
%	PARAMETERS
%
\hsize=5in
\baselineskip=12pt
\vsize=20cm
\parindent=10pt
\pretolerance=40
\predisplaypenalty=0
\displaywidowpenalty=0
\finalhyphendemerits=0
\hfuzz=2pt
\frenchspacing
\footline={\ifnum\pageno=1\else\hfil\tenrm\number\pageno\hfil\fi}
%
%	Additional Fonts
%
\input amssym.def
\def\titlefonts{\baselineskip=1.44\baselineskip
	\font\titlef=cmbx12
	\titlef
	}
\font\ninerm=cmr9
\font\ninebf=cmbx9
\font\ninei=cmmi9
\skewchar\ninei='177
\font\nineit=cmti9
\def\reffonts{\baselineskip=0.9\baselineskip
  \def\rm{\ninerm}%
	\textfont1=\ninei
  \def\bf{\ninebf}%
	\def\it{\nineit}%
	}
%
%	Formatting
%
\def\frontmatter{\vbox{}\vskip1cm\bgroup
	\leftskip=0pt plus1fil\rightskip=0pt plus1fil
	\parindent=0pt
	\parfillskip=0pt
	\pretolerance=10000
	}
\def\endfrontmatter{\egroup\bigskip}
\def\title#1{{\titlefonts#1\par}}
\def\author#1{\bigskip#1\par}
\def\address#1{\medskip{\reffonts\it#1}}
\def\email#1{\smallskip{\reffonts{\it E-mail: }\rm#1}}
\def\thanks#1{\footnote{}{\reffonts\rm\noindent#1\hfil}}
\def\section#1\par{\ifdim\lastskip<\bigskipamount\removelastskip\fi
	\penalty-250\bigskip
	\vbox{\leftskip=0pt plus1fil\rightskip=0pt plus1fil
	\parindent=0pt
	\parfillskip=0pt
  \pretolerance=10000{\bf#1}}\nobreak\medskip
	}
\def\emph#1{{\it#1}\/}
\def\proclaim#1. {\medbreak\bgroup{\noindent\bf#1.}\ \it}
\def\endproclaim{\egroup
	\ifdim\lastskip<\medskipamount\removelastskip\medskip\fi}
\newdimen\itemsize
\def\setitemsize#1 {{\setbox0\hbox{#1\ }
	\global\itemsize=\wd0}}
\def\item#1 #2\par{\ifdim\lastskip<\smallskipamount\removelastskip\smallskip\fi
	{\leftskip=\itemsize
	\noindent\hskip-\leftskip
	\hbox to\leftskip{\hfil\rm#1\ }#2\par}\smallskip}
\def\Proof#1. {\ifdim\lastskip<\medskipamount\removelastskip\medskip\fi
	{\noindent\it Proof\if\space#1\space\else\ \fi#1.}\ }
\def\endproof{\hfill\hbox{}\quad\hbox{}\hfill\llap{$\square$}\medskip}
\def\Remark. {\ifdim\lastskip<\medskipamount\removelastskip\medskip\fi
        {\noindent\bf Remark. }}

%
%		Citations
%
\newcount\citation
\newtoks\citetoks
\def\citedef#1\endcitedef{\citetoks={#1\endcitedef}}
\def\endcitedef#1\endcitedef{}
\def\citenum#1{\citation=0\def\curcite{#1}%
	\expandafter\checkendcite\the\citetoks}
\def\checkendcite#1{\ifx\endcitedef#1?\else
	\expandafter\lookcite\expandafter#1\fi}
\def\lookcite#1 {\advance\citation by1\def\auxcite{#1}%
	\ifx\auxcite\curcite\the\citation\expandafter\endcitedef\else
	\expandafter\checkendcite\fi}
\def\cite#1{\makecite#1,\cite}
\def\makecite#1,#2{[\citenum{#1}\ifx\cite#2]\else\expandafter\clearcite\expandafter#2\fi}
\def\clearcite#1,\cite{, #1]}
%
%		Reference section
%
\def\references{\section References\par
	\bgroup
	\parindent=0pt
	\reffonts
	\rm
	\frenchspacing
	\setbox0\hbox{99. }\leftskip=\wd0
	}
\def\endreferences{\egroup}
\newtoks\authtoks
\newif\iffirstauth
\def\checkendauth#1{\ifx\auth#1%
    \iffirstauth\the\authtoks
    \else{} and \the\authtoks\fi,%
  \else\iffirstauth\the\authtoks\firstauthfalse
    \else, \the\authtoks\fi
    \expandafter\nextauth\expandafter#1\fi
}
\def\nextauth#1,#2;{\authtoks={#1 #2}\checkendauth}
\def\auth#1{\nextauth#1;\auth}
\newif\ifinbook
\newif\ifbookref
\def\nextref#1 {\par\hskip-\leftskip
	\hbox to\leftskip{\hfil\citenum{#1}.\ }%
	\initnextref}
\def\initnextref{\bookreffalse\inbookfalse\firstauthtrue\ignorespaces}
\def\paper#1{{\it#1},}
\def\InBook#1{\inbooktrue in ``#1",}
\def\book#1{\bookreftrue{\it#1},}
\def\journal#1{#1\ifinbook,\fi}

\def\Vol#1{{\bf#1}}

\def\publisher#1{#1,}
\def\Year#1{\ifbookref #1.\else\ifinbook #1,\else(#1)\fi\fi}
\def\Pages#1{\makepages#1.}
\long\def\makepages#1-#2.#3{\ifinbook pp. \fi#1--#2\ifx\par#3.\fi#3}
%
%	Math Stuff
%
\newsymbol\square 1003
\newsymbol\varnothing 203F
\let\hrar\hookrightarrow
\let\rhu\rightharpoonup
\let\ot\otimes
\let\sbs\subset
\def\Ann{\mathop{\rm Ann}\nolimits}
\def\co#1{^{\mkern1mu{\rm co}\mkern1mu#1}}
\def\cop{^{\rm cop}}
\def\End{\mathop{\rm End}\nolimits}
\def\Hom{\mathop{\rm Hom}\nolimits}
\def\id{{\rm id}}
\def\Ker{\mathop{\rm Ker}}
\def\lann{\mathop{\rm lann}\nolimits}
\def\limdir{\mathop{\vtop{\offinterlineskip\halign{##\hskip0pt\cr\rm lim\cr
	\noalign{\vskip1pt}
	$\scriptstyle\mathord-\mskip-10mu plus1fil
	\mathord-\mskip-10mu plus1fil
	\mathord\rightarrow$\cr}}}}
\def\lng{\mathop{\rm length}\nolimits}
\def\Max{\mathop{\rm Max}\nolimits}
\def\op{^{\rm op}}
\def\opcop{^{\rm op,cop}}
\def\rann{\mathop{\rm rann}\nolimits}
\def\Sing{\mathop{\rm Sing}\nolimits}
\def\udim{\mathop{\rm udim}\nolimits}
\let\al\alpha
\let\be\beta
\let\ep\varepsilon
\let\ph\varphi
\let\De\Delta
\def\0{_{(0)}}
\def\1{_{(1)}}
\def\2{_{(2)}}
\def\3{_{(3)}}
\def\4{_{(4)}}
\def\C{{\cal C}}
\def\E{{\cal E}}
\def\F{{\cal F}}
\def\G{{\cal G}}
\def\Hd{H^\circ}
\def\Sd{S^\circ}

\citedef
An92
Coh86
Goo
Goo-W
Ma91
Ma92
Ma-W94
Mc-R
Mo
Mo-Sch95
Nich-Z89
Rad77
Scha00
Schn92
Schn93
Sk06
Sk07
Sk08
Sk10
Sk11
Sk-Oy06
St
Tak79
Wu-Zh03
\endcitedef

\frontmatter

\title{Flatness of Noetherian Hopf algebras\break over coideal subalgebras}
\author{Serge Skryabin}
\address{Institute of Mathematics and Mechanics,
Kazan Federal University,\break
Kremlevskaya St.~18, 420008 Kazan, Russia}
\email{Serge.Skryabin@kpfu.ru}
%\thanks{Supported by the RFBR grant 10-01-00431}

\endfrontmatter

\section
Introduction

All algebras will be considered over a fixed field $k$. The structure of Hopf 
algebras as modules over Hopf subalgebras and, more generally, over coideal 
subalgebras is of fundamental importance. Freeness results on the module 
structure exist for pointed Hopf algebras \cite{Ma91}, \cite{Rad77} and for 
finite dimensional ones \cite{Ma92}, \cite{Nich-Z89}, \cite{Sk07}. But the 
property of being a free module turns out to be too strong for other classes 
of Hopf algebras.

Commutative Hopf algebras are projective generators as modules over Hopf 
subalgebras (Takeuchi \cite{Tak79}) and are flat over right coideal subalgebras 
(Masuoka and Wigner \cite{Ma-W94}). Building upon the ideas of these papers 
Schneider \cite{Schn93} proved that any left or right Noetherian Hopf algebra 
is a faithfully flat module over central Hopf subalgebras. Some conditions on 
the algebras in conclusions of this kind are inevitable. Schauenburg 
\cite{Scha00} gave examples of Hopf algebras which are not faithfully flat 
over some Hopf subalgebras.

An algebra is said to be \emph{residually finite dimensional} if its ideals 
of finite codimension have zero intersection. Many classes of Hopf algebras 
satisfy this condition. Any residually finite dimensional Hopf algebra is flat 
over central right coideal subalgebras, and there are considerably better 
results in the case of Hopf subalgebras (see \cite{Sk08}). This shows once 
again that dealing with coideal subalgebras incurs extra complications.

The restriction to central subalgebras is clearly a serious limitation when it 
comes to noncommutative Hopf algebras. Unfortunately, the technique of central 
localizations used in \cite{Sk08} is not applicable in other situations.
The main result of the present paper is

\proclaim
Theorem 4.5.
Let $A$ be a right Noetherian right coideal subalgebra of a residually finite 
dimensional Noetherian Hopf algebra $H$. Then $A$ has a right Artinian classical 
right quotient ring{\rm,} and $H$ is left $A$-flat. Moreover{\rm,} if $A$ is a 
Hopf subalgebra{\rm,} then $H$ is left and right faithfully $A$-flat.
\endproclaim

The relevance of the classical quotient rings (the Ore rings of fractions) to 
the question of flatness has been made clear in another article \cite{Sk10}. 
We will recall that result in Theorem 4.4 of the present paper, providing a 
more direct proof of the desired conclusion. It shows immediately that Theorem 
4.5 holds when $A$ and $H$ are additionally assumed to be semiprime since then 
the classical quotient rings of $A$ and $H$ exist by the Goldie Theorem. 
Without this additional condition it is not easy to establish the existence of 
the classical quotient rings, and here lies the main problem since the usual 
methods do not work.

The right coideal subalgebras of a Hopf algebra $H$ are module algebras for 
the dual Hopf algebra $\Hd$. This suggests a reformulation of the problem in 
terms of module algebras. Inasmuch as the quotient rings are concerned, 
switching to the module structures is essential since those extend to the 
quotient rings, while the comodule structures generally do not.

Let now $A$ be a right Noetherian $H$-semiprime $H$-module algebra. The first 
attempt to deal with its quotient ring was not fully successful. In 
\cite{Sk-Oy06} it was shown that $A$ has a semiprimary generalized quotient 
ring $Q$ constructed with respect to a certain filter of right ideals. The 
property of being semiprimary is close to being Artinian, but still it does 
not seem to allow one to deduce that $Q$ is a classical quotient ring. The 
latter conclusion was obtained in \cite{Sk-Oy06} only for some classes of Hopf 
algebras. In the present paper we will prove it assuming that the action of 
$H$ on $A$ is \emph{locally finite}, i.e., each element of $A$ is contained in 
a finite dimensional $H$-submodule:

\proclaim
Theorem 4.1.
Let $A$ be a right Noetherian $H$-semiprime $H$-module algebra such that the 
action of $H$ on $A$ is locally finite. Then $A$ has a right Artinian 
classical right quotient ring.
\endproclaim

This result is sufficient to derive Theorem 4.5 since the action of $\Hd$ on 
$H$, and therefore on all right coideal subalgebras of $H$, is locally finite.

It should be stressed that there are no restrictions on the Hopf algebra $H$ 
in Theorem 4.1. To achieve this generality we have to revise the approach of 
\cite{Sk-Oy06} where the antipode of $H$ was assumed to be bijective. Using a 
slightly modified filter $\E'_H$ of right ideals, we are still able to prove 
that the corresponding quotient ring $Q$ is semiprimary and $H$-semiprime. 
This is done in the first two sections of the present paper.

However, we do not need other parts of \cite{Sk-Oy06} since we provide 
completely different arguments to analyze the structure of $Q$ in section 3 of 
the paper. In particular, we needn't bother with the selfinjectivity of $Q$. 
The local finiteness of the action leads very quickly to the decomposition of 
$Q$ as a direct product of $H$-simple algebras. Then we show that each 
$Q$-module has no nonzero $\E'_H$-torsion elements, which is a crucial 
property in the verification that $Q$ is indeed a classical quotient ring.

The final results are presented in section 4 of the paper. Most of them have 
been discussed already in this introduction. Combining our approach here with 
an already known result on the antipode proved in \cite{Sk06} we also obtain

\proclaim
Theorem 4.3.
Let $H$ be either right or left Noetherian residually finite dimensional Hopf 
algebra. Then its antipode $S:H\to H$ is bijective. Hence $H$ is right and 
left Noetherian simultaneously.
\endproclaim

By Theorem 4.5 applied to $A=H$ each residually finite dimensional Noetherian 
Hopf algebra $H$ has an Artinian classical quotient ring. In the case when $H$ 
is a Noetherian affine PI Hopf algebra such a conclusion was deduced earlier 
by Wu and Zhang \cite{Wu-Zh03} as a consequence of Gorensteinness of $H$. As a 
matter of fact, the assumption of Theorem 4.5 is satisfied in this case, and 
so we obtain an alternative proof. Indeed, it was proved by Anan'in 
\cite{An92} that each right Noetherian finitely generated PI algebra is 
residually finite dimensional. One may wonder whether every Noetherian Hopf 
algebra is necessarily residually finite dimensional.

\bigbreak
\centerline{\it Terminology and Notation}

\medskip
For a subset $X$ of a ring $R$ we denote by $\lann_RX$ and $\rann_RX$, 
respectively, the left and right annihilators of $X$ in $R$. An element $s\in R$ 
is called \emph{right regular} if $\rann_Rs=0$. Left regular elements are 
defined by the condition $\lann_Rs=0$, and $s$ is called \emph{regular} if it 
is both right and left regular.

A ring $Q$ containing $R$ as a subring is said to be a \emph{classical right 
quotient ring} of $R$ if all regular elements of $R$ are invertible in $Q$ and 
each element of $Q$ can be written as $as^{-1}$ for some $a,s\in R$ with $s$ 
being regular. See \cite{Goo-W} or \cite{Mc-R} for information on related 
topics.

A ring is called \emph{semiprimary} if its Jacobson radical is nilpotent 
and the factor ring by the Jacobson radical is semisimple Artinian.

For general facts and definitions concerning Hopf algebras we refer to 
\cite{Mo}. Let $H$ be a Hopf algebra over $k$. We denote by $\De$, $\ep$, $S$ 
its comultiplication, counit and antipode, and we write $\De(h)=\sum h\1\ot 
h\2\in H\ot H$ for $h\in H$.

A \emph{right coideal} of $H$ is any subspace $U$ such that $\De(U)\sbs U\ot 
H$. A subalgebra of $H$ satisfying this condition is called a \emph{right 
coideal subalgebra}.

All algebras are assumed to be associative and unital. An \emph{$H$-module 
algebra} $A$ is equipped with a left $H$-module structure such that
$$
h(ab)=\sum\,(h\1a)\,(h\2b)\quad\hbox{for all $h\in H$ and $a,b\in A$}.
$$
The following two useful identities hold in such an algebra:
$$
(ha)b=\sum\,h\1\bigl(a\,S(h\2)b\bigr),\qquad
a\bigl(S(h)b\bigr)=\sum\,S(h\1)\bigl((h\2a)\,b\bigr).
$$
It follows that $\lann_AV$ is an $H$-submodule of $A$ for each $S(H)$-submodule 
$V$. If $V$ is an $H$-submodule, then $\rann_AV$ is an $S(H)$-submodule, but 
we cannot be sure that $\rann_AV$ is an $H$-submodule unless $S(H)=H$ (cf. 
\cite{Coh86, Cor. 2}). Similarly, the left annihilators of $S^2(H)$-submodules 
of $A$ are $S(H)$-submodules.

An $H$-module algebra $A$ is called \emph{$H$-simple} if $A\ne0$ and $A$ has 
no $H$-stable ideals except the zero ideal and the whole $A$.
An $H$-module algebra $A$ is \emph{$H$-prime} if $A\ne0$ and $IJ\ne0$ for all 
nonzero $H$-stable ideals $I$ and $J$ of $A$. And $A$ is \emph{$H$-semiprime} 
if $A$ contains no nonzero nilpotent $H$-stable ideals. By an ideal we mean a 
two-sided ideal. The action of $H$ on $A$ is said to be \emph{locally finite} 
if $\,\dim Ha<\infty\,$ for all $a\in A$.

\section
1. The filter of right ideals

Recall that a (right) {\it Gabriel topology} on a ring $R$ is any set $\G$ of 
right ideals of $R$ satisfying the four conditions listed below where $I,\,J$ 
are assumed to be right ideals of $R$ and we use the notation $(I:a)=\{x\in 
R\mid ax\in I\}$:

\setitemsize(T4)
\item(T1)
If $J\in\G$ and $J\sbs I$ then $I\in\G${\rm;}

\item(T2)
If $I,J\in\G$ then $I\cap J\in\G${\rm;}

\item(T3)
If $I\in\G$ then $(I:a)\in\G$ for each $a\in R${\rm;}

\item(T4)
If $J\in\G$ and $(I:a)\in\G$ for all $a\in J$ then $I\in\G$.

\smallskip
With a Gabriel topology $\G$ one associates a hereditary torsion theory (see 
\cite{St, Ch. VI, Th. 5.1}). A right $R$-module is said to be \emph{$\G$-torsion} 
if each of its elements is annihilated by a right ideal in $\G$. The class of 
$\G$-torsion modules is closed under submodules, factor modules, coproducts, 
and extensions. An arbitrary right $R$-module $V$ has a largest $\G$-torsion 
submodule. This submodule consists of all elements of $V$ whose annihilators 
in $R$ belong to $\G$. A right $R$-module is called \emph{$\G$-torsionfree} if 
it contains no nonzero $\G$-torsion submodules.

Let $A$ be a left $H$-module algebra. Denote by $\E$ the set of all essential 
right ideals of $A$. Recall that a right ideal is said to be \emph{essential} 
if it has nonzero intersection with each nonzero right ideal. It is well-known 
that $\E$ satisfies (T1)--(T3). 

In \cite{Sk-Oy06} we worked with the set $\E_H$ of all right ideals $I$ of $A$ 
such that for each $h\in H$ one has $hJ\sbs I$ for some $J\in\E$. However, in 
the case when $S(H)\ne H$ we do not get the necessary properties of this 
filter. For this reason we will use a slightly different filter of right 
ideals. Note that $S(H)$ is a Hopf subalgebra of $H$ since the antipode 
$S:H\to H$ is a Hopf algebra antiendomorphism.

Denote by $\E'_H$ the set of right ideals $I$ of $A$ having the property that 
for each $h\in S(H)$ one has $hJ\sbs I$ for some right ideal $J\in\E$ 
depending on $I$ and $h$. We will write $\E'_H(A)$ instead of $\E'_H$ when we 
need to indicate the algebra $A$.

Since $1\in S(H)$, each right ideal $I\in\E'_H$ contains an essential right 
ideal, and therefore is itself essential. So $\E_H\sbs\E'_H\sbs\E$. Clearly 
$\E'_H=\E_H$ when $S$ is surjective.

\medskip
For a coalgebra $C$ denote by $[C,A]$ the vector space $\Hom_k(C,A)$ equipped 
with the convolution multiplication. If $\dim C<\infty$, then $[C,A]\cong A\ot 
C^*$ as algebras, and if $C\sbs H$, there is an algebra homomorphism 
$\tau:A\to[C,A]$ defined by the rule $\tau(a)(c)=ca$. One can check that 
$[C,A]=\tau(A)C^*$, and so $[C,A]$ is finitely generated as a left 
$\tau(A)$-module. When $S$ is not bijective, we cannot derive the right hand 
version of this conclusion. In order to use the finiteness property in one of 
the later arguments, we have to modify the previous construction. 

Denote by $\F$ the set of all finite dimensional subcoalgebras of $H$. 
Let $C\in\F$, and let $C\cop$ be $C$ with the opposite comultiplication. The 
algebra $[C\cop,A]$ is defined on the vector space $\Hom_k(C,A)$ by means of 
the multiplication
$$
(\xi\times\eta)(c)=\sum\xi(c\2)\eta(c\1),\qquad  
\xi,\eta\in\Hom_k(C,A),\quad c\in C.
$$
Clearly $[C\cop,A]\cong A\ot(C^*)\op$. Define a map $\rho:A\to[C\cop,A]$, 
$\,a\mapsto\rho_a$, setting
$$
\rho_a(c)=S(c)a,\qquad a\in A,\ c\in C.
$$
This map is an algebra homomorphism since
$$
\rho_{ab}(c)=S(c)(ab)=\sum\,\bigl(S(c\2)a\bigr)\bigl(S(c\1)b\bigr)
=(\rho_a\times\rho_b)(c)
$$
for all $a,b\in A$ and $c\in C$.

\proclaim
Lemma 1.1.
For $C\in\F$ the algebra $[C\cop,A]$ is a free $A$-module of finite rank with 
respect to the right action of $A$ obtained via $\rho:$
$$
\xi\cdot_\rho a=\xi\times\rho_a\quad\hbox{where $\,\xi\in\Hom_k(C,A),$ $\,a\in A$}.
$$
In particular{\rm,} $\rho$ is injective{\rm,} and so the subalgebra 
$\rho(A)\sbs[C\cop,A]$ is isomorphic to $A${\rm,} whenever $C\ne0$.
\endproclaim

\Proof.
Clearly $\Hom_k(C,A)\cong C^*\ot A$ is a free $A$-module of finite rank with 
respect to another right action of $A$ such that
$$
(\xi a)(c)=\xi(c)a\quad\hbox{for $\xi\in\Hom_k(C,A)$, $\,a\in A$, $\,c\in C$}.
$$
So it suffices to check that $\,\cdot_\rho\,$ is an isomorphic $A$-module 
structure. Define a linear transformation $\Phi$ of $\Hom_k(C,A)$ setting
$$ 
(\Phi\xi)(c)=\sum S(c\1)\xi(c\2).
$$
Since
$$
\eqalign{
\bigl(\Phi(\xi a)\bigr)(c)
=\sum S(c\1)\bigl(\xi(c\2)a\bigr)
&{}=\sum\,\bigl(S(c\2)\xi(c\3)\bigr)\,\bigl(S(c\1)a\bigr)\cr
&\qquad{}=\sum\,(\Phi\xi)(c\2)\,\rho_a(c\1)=(\Phi\xi\times\rho_a)(c)
}
$$
for all $c\in C$, we get $\Phi(\xi a)=\Phi(\xi)\cdot_\rho a$ for all 
$\xi\in\Hom_k(C,A)$ and $a\in A$. The inverse transformation $\Phi^{-1}$ is 
defined by the rule
$$
(\Phi^{-1}\xi)(c)=\sum c\1\xi(c\2).
$$
Thus $\Phi$ is bijective, and so $\Phi$ is indeed an isomorphism between the 
two $A$-module structures on $\Hom_k(C,A)$.
\endproof

For any right ideal $I$ of $A$ and a subcoalgebra $C\sbs H$ put
$$
I_C=\tau^{-1}\bigl(\Hom_k(C,I)\bigr)=\{x\in A\mid Cx\sbs I\}.
$$
Since $\tau:A\to[C,A]$ is an algebra homomorphism and $\Hom_k(C,I)$ is a right 
ideal of $[C,A]$, it is clear that $I_C$ is a right ideal of $A$. Note that
$$
I_{S(C)}=\rho^{-1}\bigl(\Hom_k(C,I)\bigr)=\{x\in A\mid S(C)x\sbs I\}.
$$

\proclaim
Lemma 1.2.
A right ideal $I$ of $A$ is in $\E'_H$ if and only if $I_{S(C)}\in\E$ for each
$C\in\F$. Moreover{\rm,} $I_{S(C)}\in\E'_H$ whenever $I\in\E'_H$.
\endproclaim

\Proof.
Suppose that $I\in\E'_H$. Given $C\in\F$ and $h\in S(H)$, let $X$ be any basis 
of the finite dimensional subspace $S(C)h\sbs S(H)$. Since $X$ is finite and 
$\E$ is closed under finite intersections of right ideals, there exists 
$J\in\E$ such that $gJ\sbs I$ for all $g\in X$. Then $S(C)hJ\sbs I$, that is, 
$hJ\sbs I_{S(C)}$. This establishes the inclusion $I_{S(C)}\in\E'_H\sbs\E$. 

Conversely, since every element of $H$ is contained in a finite dimensional 
subcoalgebra, $S(H)$ is the union of the subcoalgebras $S(C)$ with $C\in\F$, 
and obviously we have $hI_{S(C)}\sbs I$ for all $h\in S(C)$. This shows that 
$I\in\E'_H$ whenever $I_{S(C)}\in\E$ for all $C\in\F$.
\endproof

\proclaim
Lemma 1.3.
The set $\E'_H$ satisfies {\rm(T1)--(T3)}.
\endproclaim

\Proof.
Properties (T1), (T2) for $\E'_H$ follow easily from the respective 
properties of $\E$ since $(I\cap J)_C=I_C\cap J_C$ and, in particular, 
$I_C\sbs J_C$ whenever $I\sbs J$.

Let us check (T3). Let $I\in\E'_H$ and $a\in A$. For each $C\in\F$ the right 
ideal $I_{S(C)}$ is essential by Lemma 1.2, and we have to show that so is
$$
(I:a)_{S(C)}=\{x\in A\mid a\,\bigl(S(C)x\bigr)\sbs I\}.
$$
Put $\,K=\{x\in A\mid(Ca)x\sbs I_{S(C)}\}\,$.
Then $K\in\E$ since $Ca$ is a finite dimensional subspace of $A$, and $\E$ 
satisfies (T2), (T3). To complete the proof it remains to show that 
$K\sbs (I:a)_{S(C)}$. This containment does hold because
$$
a\,\bigl(S(c)x\bigr)=\sum S(c\1)\bigl((c\2a)x\bigr)\in S(C)\,I_{S(C)}\sbs I
$$
for all $x\in K$ and $c\in C$.
\endproof

\proclaim
Proposition 1.4.
Suppose that $A$ is $H$-semiprime and satisfies ACC on right annihilators. 
Then $\E'_H$ is a Gabriel topology and $A$ is $\E'_H$-torsionfree as a module 
over itself with respect to right multiplications.
\endproclaim

\Proof.
Since $\E'_H$ satisfies (T1)--(T3), the set
$$
N=\{a\in A\mid\rann_A a\in\E'_H\}
$$
is a right ideal of $A$. Obviously, $N$ is stable also under left 
multiplications in $A$. Suppose $a\in N$. Then $aI=0$ for some $I\in\E'_H$. If 
$C\in\F$, then $I_{S(C)}\in\E'_H$ by Lemma 1.2. Now
$$
(ca)x=\sum c\1\bigl(a\,S(c\2)x\bigr)=0
$$
for all $c\in C$ and $x\in I_{S(C)}$ since $S(C)x\sbs I$. Thus $bI_{S(C)}=0$ 
for each $b\in Ca$.  This shows that $Ca\sbs N$. Since each element of $H$ is 
contained in a finite dimensional subcoalgebra, we conclude that $N$ is an 
$H$-stable two-sided ideal of $A$.

Recall that the right singular ideal $\Sing A$ of $A$ is a two-sided ideal 
consisting of all elements of $A$ whose right annihilators are essential right 
ideals. Since $\E'_H\sbs\E$, we have $N\sbs\Sing A$. According to \cite{Mc-R, 
Lemma 2.3.4} the ascending chain condition on right annihilators implies that 
$\Sing A$ is nilpotent. Hence so too is $N$, and the $H$-semiprimeness of $A$ 
yields $N=0$.

Vanishing of $N$ means that $\lann_AI=0$ for each $I\in\E'_H$. In other 
words, $A$ is $\E'_H$-torsionfree. It remains to verify that $\E'_H$ satisfies 
(T4).

Let $I$ and $J$ be two right ideals of $A$ such that $J\in\E'_H$ and 
$(I:a)\in\E'_H$ for all $a\in J$. We have to show that $I\in\E'_H$. Let 
$C\in\F$. Then $J_{S(C)}\in\E$, and we will check that $I_{S(C)}\in\E$ too. 
For this we have to show that $I_{S(C)}\cap R\ne0$ for each nonzero right 
ideal $R$ of $A$. But $J_{S(C)}\cap R\ne0$, so that it suffices to consider 
the right ideals of the form $R=bA$ where $0\ne b\in J_{S(C)}$. Fix such 
an element $b$ and put
$$
K=\{x\in A\mid\bigl(S(C)b\bigr)\,x\sbs I\}.
$$
Here $S(C)b$ is a finite dimensional subspace of $J$. Taking its basis, say 
$b_1,\ldots,b_n$, we get $K=\bigcap_{i=1}^n(I:b_i)\in\E'_H$ since 
$b_1,\ldots,b_n\in J$. If $y\in K_{S(C)}$, then $S(C)y\sbs K$, and therefore
$$
S(c)(by)=\sum\bigl(S(c\2)b\bigr)\,\bigl(S(c\1)y\bigr)\sbs
\bigl(S(C)b\bigr)\,K\sbs I
$$
for all $c\in C$, i.e. $by\in I_{S(C)}$. We see that
$$
b\,K_{S(C)}\sbs I_{S(C)}\cap bA.
$$
But $K_{S(C)}\in\E'_H$ by Lemma 1.2. As we have proved already, all right 
ideals in $\E'_H$ have zero left annihilators. Hence $\,bK_{S(C)}\ne0$, and 
therefore $I_{S(C)}\cap bA\ne0$.
\endproof

Later we will have to work with $H$-module algebras which are not right 
Noetherian, but only right Goldie. The ACC on right annihilators is one of 
Goldie conditions. The second one is the ACC on direct sums of right ideals, 
which can be interpreted as the finiteness of the right uniform dimension. Our 
next aim is to show that in the presence of the Goldie conditions the filter 
$\E'_H$ is sufficiently large (see Lemma 1.6).

Recall that the \emph{uniform dimension} $\,\udim M\,$ of a module $M$ is the 
largest number of nonzero submodules forming a direct sum, and $\,\udim 
M<\infty\,$ if no infinite direct sum of nonzero submodules exist. If $R$ is a 
subring of a ring $T$, let $T_R$ be $T$ regarded as a right 
$R$-module with respect to the action of $R$ on $T$ by right multiplications.  
Our argument is based on the following ring-theoretic observation:

\proclaim
Lemma 1.5.
Suppose that $R$ is a subring of a ring $T$ such that $\udim T_R<\infty$. Then
$xT+\rann_Tx$ is an essential submodule of $T_R$ for any element $x\in T$ with
$\rann_Tx=\rann_Tx^2$.
\endproclaim

\Proof.
Suppose that $V\sbs T_R$ is a submodule such that $V\cap(xT+\rann_Tx)=0$. Then 
$x^iV\cong V$ for each $i\ge0$, and the sum $\sum_{i=0}^\infty x^iV\sbs T_R$ 
is direct. The finiteness of the uniform dimension entails $V=0$. 
\endproof

\proclaim
Lemma 1.6.
Suppose that $\udim A_A<\infty$. Then
$$
uA+\rann_Au\in\E'_H
$$
for any element $u\in A$ with $\rann_Au=\rann_Au^2$. In particular{\rm,} 
$uA\in\E'_H$ whenever $u$ is right regular in $A$.
\endproclaim

\Proof.
Put $I=uA+\rann_Au$. We want to apply Lemma 1.5 with $T=[C\cop,A]$ and 
$R=\rho(A)$ where $C\in\F$. By Lemma 1.1 $R\cong A$ and $T$ is a free right 
$R$-module of finite rank. Hence $\udim T_R<\infty$.

Making use of the identification $T\cong A\ot(C^*)\op\,$ let $x=u\ot1\in T$. 
We have $\,xT\cong uA\ot(C^*)\op\,$ and 
$\,\rann_Tx\cong(\rann_Au)\ot(C^*)\op$. Hence
$$
xT+\rann_Tx=\Hom_k(C,I)\cong I\ot(C^*)\op.
$$
Since $x^2=u^2\ot1$, we deduce that 
$\,\rann_Tx^2=(\rann_Au^2)\ot(C^*)\op=\rann_Tx$. All assumptions of Lemma 1.5 
thus hold, and $\Hom_k(C,I)$ is then an essential right $R$-submodule of $T$.

Recall that
$$
I_{S(C)}=\rho^{-1}\bigl(\Hom_k(C,I)\bigr).
$$
Let $0\ne b\in A$. If $\rho_b=0$ in $T$ then $b\in I_{S(C)}$ 
since $\rho_b\in\Hom_k(C,I)$. If $\rho_b\ne0$ then
$$
\rho_bR\cap\Hom_k(C,I)\ne0;
$$
since $\rho_bR=\rho(bA)$, there exists $a\in bA$ such that 
$0\ne\rho_a\in\Hom_k(C,I)$. In the latter case $0\ne a\in I_{S(C)}$. Thus 
$I_{S(C)}\cap bA\ne0$ in any case, and so $I_{S(C)}\in\E$. Lemma 1.2 
completes the proof.
\endproof

As we have seen, under the hypothesis of Proposition 1.4 all right ideals in 
$\E'_H$ have zero left annihilators. It will be important later that the right 
annihilators are zero as well, but we can prove this only under stronger 
restrictions:

\proclaim
Lemma 1.7.
Suppose that $A$ is $S^2(H)$-semiprime and satisfies ACC on right 
annihilators. Then\/ $\rann_AI=0$ for each $I\in\E'_H$.
\endproclaim

\Proof.
For each right ideal of $A$ its right annihilator in $A$ is a two-sided ideal. 
By the ACC the set
$$
\{\rann_AI\mid I\in\E'_H\}
$$
has a maximal element, say $K$. But this set is directed by inclusion since 
the set $\E'_H$ is directed by inverse inclusion according to property (T2) 
and since the correspondence $I\mapsto\rann_AI$ reverses inclusions. Therefore 
$K$ is the largest among all right annihilators of right ideals in $\E'_H$. 
We have to show that $K=0$.

Now pick $I\in\E'_H$ such that $K=\rann_AI$. If $a\in K$, $C\in\F$, $c\in C$, 
$x\in I_{S(C)}$, then
$$
x\,\bigl(S^2(c)a\bigr)=\sum\,S^2(c\2)\bigl((S(c\1)x)\,a\bigr)=0
$$
since $S(C)x\sbs I$ and $Ia=0$. This shows that $S^2(c)a\in\rann_AI_{S(C)}$. 
Lemma 1.2 tells us that $I_{S(C)}\in\E'_H$. Hence $\rann_AI_{S(C)}\sbs K$ by 
the choice of $K$, and it follows that $S^2(C)K\sbs K$. Since $H$ is the union 
of subcoalgebras $C\in\F$, we conclude that $K$ is stable under the action of 
$S^2(H)$.

The left annihilator $L=\lann_AK$ is a two-sided ideal as well. It is stable 
under the action of $S(H)$ since $K$ is an $S^2(H)$-submodule of $A$. Hence 
$K\cap L$ is an $S^2(H)$-stable ideal. Since $(K\cap L)^2\sbs LK=0$, we deduce 
that $K\cap L=0$.

Since $KL\sbs K\cap L$, it follows that $KL=0$ too. But $L\in\E'_H$ since $L$ 
contains any right ideal $I\in\E'_H$ such that $K=\rann_AI$. By Proposition 
1.4 $\,\lann_AL=0$, which entails $K=0$, as required.
\endproof

\proclaim
Lemma 1.8.
Let $I$ be an ideal of finite codimension in $H$. Then $S^n(H)+I=H$ for all 
$n>0$.
\endproclaim

\Proof.
Recall that the finite dual $\Hd$ of $H$ is a Hopf algebra consisting of all 
linear functions $H\to k$ vanishing on an ideal of finite codimension in $H$. 
The antipode $\Sd$ of $\Hd$ is defined by the rule $\Sd(f)=f\circ S$ for 
$f\in\Hd$.

By \cite{Sk06, Th. A} $\Sd$ is injective since $\Hd$ is always residually finite 
dimensional and, as a consequence, weakly finite. If $f\in H^*$ is any linear 
function vanishing on $S^n(H)+I$, then $f\in\Hd$ and 
${\Sd}^{\mskip0.5mu n}(f)=0$, 
whence $f=0$ by injectivity of $\Sd$. This yields the desired conclusion.
\endproof

\proclaim
Corollary 1.9.
Let $V$ be a locally finite dimensional $H$-module{\rm,} so that each element 
of $V$ is contained in a finite dimensional submodule. Then each 
$S^n(H)$-submodule of $V$ is an $H$-submodule.
\endproclaim

\Proof.
Let $I$ be the annihilator of the $H$-submodule $Hv$ generated by some element 
$v\in V$. Since $\dim Hv<\infty$, this ideal of $H$ has finite codimension. By 
Lemma 1.8 $S^n(H)+I=H$, and we get $Hv=S^n(H)v$ since $Iv=0$. Hence $Hv$ is 
contained in each $S^n(H)$-submodule of $V$ containing $v$.
\endproof
\proclaim
Corollary 1.10.
Suppose that $A$ is $H$-semiprime and satisfies ACC on right annihilators. 
If the action of $H$ on $A$ is locally finite{\rm,} then\/ $\rann_AI=0$ 
for each $I\in\E'_H$.
\endproclaim

\Proof.
By Corollary 1.9 each $S^2(H)$-stable ideal of $A$ is stable under the action 
of the whole $H$. Therefore there is no difference between the 
$H$-semiprimeness and the $S^2(H)$-semiprimeness of $A$, and so Lemma 1.7 
applies.
\endproof

\section
2. The quotient ring

Let $R$ be a ring and $\G$ a filter of right ideals satisfying the axioms 
(T1)--(T4) of a Gabriel topology. The abelian groups $\Hom_R(I,R)$ with 
$I\in\G$ form an inductive system, and in the case when $R$ is $\G$-torsionfree 
as a right $R$-module the {\it localization} of $R$ with respect to $\G$ is 
defined as the limit
$$
R_\G=\limdir_{I\in\G}\Hom_R(I,R).
$$
If $\al:I\to R$ and $\be:J\to R$ are two right $R$-linear maps, where $I,J\in\G$, 
then $\al\circ\be$ is defined on $\be^{-1}(I)$. The $R$-module $J/\be^{-1}(I)$ 
is $\G$-torsion since it embeds in the $\G$-torsion module $R/I$. Since $R/J$ 
is $\G$-torsion, so is $R/\be^{-1}(I)$ too, which means that $\be^{-1}(I)\in\G$. 
Thus $\al\circ\be$ represents an element of $R_\G$, and this one is taken to 
be the product of the two elements represented by $\al$ and $\be$ 
respectively. In this way $R_\G$ acquires a ring structure. We call $R_\G$ 
with this structure \emph{the quotient ring} of $R$ with respect to $\G$.

The ring $R$ is identified with the subring of $R_\G$ consisting of all 
elements represented by left multiplications in $R$. If $q\in R_\G$ is an 
arbitrary element represented by $\al:I\to R$, where $I\in\G$, then 
$qx=\al(x)$ for all $x\in I$; hence $qI\sbs R$, and $qI\ne0$ unless $q=0$. In 
particular, each nonzero right $R$-submodule of $R_\G$ has a nonzero 
intersection with $R$, so that $R_\G$, regarded as a right $R$-module, is an 
essential extension of $R$. It follows that $R_\G$, along with $R$, is 
$\G$-torsionfree.

Suppose now that $\G$ is a Gabriel topology on a left $H$-module algebra $A$ 
such that $A$ is $\G$-torsionfree. The right ideals in $\G$ form a 
neighbourhood base of 0 for a topology making $A$ into a topological algebra. 
If all elements of $H$ operate on $A$ as continuous transformations, then the 
action of $H$ on $A$ is said to be $\G$-continuous, and it is known in this 
case that the action extends to the quotient ring $A_\G$ \cite{Mo-Sch95, 
Th. 3.13}. It will be important for us that the conclusion of that theorem 
remains valid under a slightly weaker assumption when continuity of the action 
is required only for elements of $S(H)$:

\proclaim
Lemma 2.1.
Suppose that all elements of $S(H)$ operate on $A$ as $\G$-continuous 
transformations. Then $A_\G$ is a left $H$-module algebra with respect to an 
action of $H$ extending the given action on $A$.
\endproclaim

\Proof.
The continuity assumption means that for each $h\in S(H)$ and each $I\in\G$ 
there exists $I_h\in\G$ such that $hI_h\sbs I$. If $C\in\F$, then $S(C)$ is a 
finite dimensional subspace of $S(H)$; therefore by (T2) for each $I\in\G$ 
there exists $K\in\G$ such that $hK\sbs I$ for all $h\in S(C)$, i.e. $K\sbs 
I_{S(C)}$ in the notation of section 1. Note that the latter inclusion and 
(T1) imply that $I_{S(C)}\in\G$.

Given any $h\in H$ and a right $A$-linear map $\al:I\to A$ where $I\in\G$, let 
$C_h\in\F$ be the smallest subcoalgebra containing $h$, and define 
$h\al:I_{S(C_h)}\to A$ by the rule
$$
(h\al)(x)=\sum h\1\,\al\bigl(S(h\2)x\bigr),\qquad x\in I_{S(C_h)}.
$$
As in \cite{Coh86, Th. 18}, one checks that the map $h\al$ is $A$-linear. 
Since $I_{S(C_h)}\in\G$, this map represents an element of $A_\G$.

Thus for $h\in H$ and $q\in A_\G$ we can define $hq$ to be the element of 
$A_\G$ represented by $h\al$ where $\al$ is any representative of $q$. If $g$ 
is a second element of $H$, then $(gh)q=g(hq)$ since the two $A$-linear maps 
$(gh)\al$ and $g(h\al)$ agree on the right ideal
$$
\{x\in A\mid S(C_h)S(C_g)x\sbs I\}=\bigl(I_{S(C_h)}\bigr)_{S(C_g)}\in\G.
$$

Let $q,q'\in A_\G$ be represented by $A$-linear maps $\al:I\to A$ and 
$\be:J\to A$ where $I,J\in\G$. For each $c\in C_h$ the map $c\al$ is defined 
on $I_{S(C_h)}\in\G$ and $c\be$ is defined on $J_{S(C_h)}\in\G$. Since 
$\{c\be\mid c\in C_h\}$ is a finite dimensional subspace of $A$-linear maps 
$J_{S(C_h)}\to A$, there exists $K\in\G$ such that $K\sbs J_{S(C_h)}$ and 
$(c\be)(K)\sbs I_{S(C_h)}$ for all $c\in C_h$. If $x\in K$ and $c\in C_h$, 
then
$$
\be\bigl(S(c)x\bigr)=\sum S(c\1)c\2\,\be\bigl(S(c\3)x\bigr)
=\sum S(c\1)\bigl((c\2\be)(x)\bigr)\in S(C_h)I_{S(C_h)},
$$
whence $\be\bigl(S(c)x\bigr)\in I$, and so $\al\circ\be$ is defined at $S(c)x$. 
We get
$$
\eqalign{
\bigl(h(\al\circ\be)\bigr)(x)=\sum h\1\,\al\Bigl(\be\bigl(S(h\2)x\bigr)\Bigr)
&{}=\sum h\1\,\al\Bigl(S(h\2)h\3\,\be\bigl(S(h\4)x\bigr)\Bigr)\cr
&{}=\sum\,(h\1\al)\bigl((h\2\be)(x)\bigr)\cr
}
$$
Thus $h(\al\circ\be)$ agrees with $\sum\,(h\1\al)\circ(h\2\be)$ on $K$. This 
means that
$$
h(qq')=\sum\,(h\1q)(h\2q'),
$$
i.e. the $H$-module structure on $A_\G$ is compatible with the multiplication.
\endproof

It was proved in Lemma 1.2 that $I_{S(C)}\in\E'_H$ for each $I\in\E'_H$ and 
each $C\in\F$. This fact shows that all elements of $S(H)$ operate on $A$ as 
$\E'_H$-continuous transformations.

\proclaim
Proposition 2.2.
Suppose that $A$ is $H$-semiprime and right Noetherian. Then the quotient ring 
$Q$ of $A$ with respect to the filter of right ideals $\E'_H$ is a semiprimary 
$H$-semiprime $H$-module algebra.
\endproclaim

\Proof.
Note that Proposition 1.4 applies to $A$ since $A$ satisfies ACC on arbitrary 
right ideals. Hence $\E'_H$ is a Gabriel topology and $A$ is 
$\E'_H$-torsionfree. By Lemma 2.1 $Q$ is a left $H$-module algebra with 
respect to an action of $H$ extending the given action on $A$.

By the general properties of the quotient rings $Q$ is an essential extension 
of $A$ in the category of right $A$-modules. In particular, each nonzero right 
ideal of $Q$ has nonzero intersection with $A$. If $I$ is a nilpotent 
$H$-stable ideal of $Q$, then $I\cap A$ is a nilpotent $H$-stable ideal of 
$A$. Then $I\cap A=0$ by $H$-semiprimeness of $A$, and we must have $I=0$. 
Therefore $Q$ is $H$-semiprime.

The less obvious part of the proof is to show that $Q$ is semiprimary. Here 
the arguments follow \cite{Sk-Oy06} with the filter $\E_H$ replaced by $\E'_H$ 
everywhere. We indicate below the main steps referring to \cite{Sk-Oy06} for 
other details.

\proclaim
Claim 1.
$Q$ is right Goldie{\rm,} i.e.{\rm,} $Q$ satisfies ACC on right annihilators 
and has finite right uniform dimension.
\endproclaim

This is a general property of the quotient rings of right Noetherian rings 
(see \cite{Sk-Oy06, Lemma 6.1}).

\proclaim
Claim 2.
For a right ideal $I$ of $Q$ one has $I\in\E'_H(Q)$ if and only if 
$I\cap A\in\E'_H(A)$.
\endproclaim

This is an analog of \cite{Sk-Oy06, Lemma 6.2}. Setting $J=I\cap A$, we have 
$$
J_{S(C)}=I_{S(C)}\cap A\qquad\hbox{for each $C\in\F$}.
$$
Therefore $I_{S(C)}$ is an essential right ideal of $Q$ if and only if 
$J_{S(C)}$ is an essential right ideal of $A$ (see \cite{Goo, Prop. 2.32(a)}). 
It remains to apply the characterization of the filters $\E'_H(A)$, $\E'_H(Q)$ 
given in Lemma 1.2.

\proclaim
Claim 3.
If $I\in\E'_H(Q)$ then each right $Q$-linear map $I\to Q$ is induced by a 
left multiplication in $Q$.
\endproclaim

It is well-known that the quotient ring $Q$ coincides with its own 
localization with respect to the filter $\G^e$ of right ideals $I$ such that 
$I\cap A\in\E'_H(A)$ \cite{St, Ch.\ X, \S2}. By Claim 2 we have $\G^e=\E'_H(Q)$, 
and therefore bijectivity of the canonical map $Q\to Q_{\G^e}$ amounts to 
Claim 3.

\proclaim
Claim 4.
For any $u\in Q$ satisfying $\,\rann_Qu=\rann_Qu^2$ there exists an idempotent
$e\in Q$ such that $u$ is an invertible element of the ring $eQe$ with unity $e$.
\endproclaim

Put $Y=\rann_Qu$. The right ideal $I=uQ+Y$ is in $\E'_H(Q)$ by Lemma 1.6, and 
the assumption about $u$ implies that the sum here is direct. The projection 
of $I$ onto $uQ$ is a right $Q$-linear map $I\to Q$. By Claim 3 it is the 
restriction to $I$ of the left multiplication by some element $e\in Q$. Then 
$eu=u$ and $eY=0$. Since $(ue-u)I=0$ and $(e^2-e)I=0$, it follows 
that $ue=u$ and $e^2=e$ by torsionfreeness of $Q$.

There is also a right $Q$-linear map $I\to Q$ such that $uq\mapsto eq$ for all 
$q\in Q$ and $q\mapsto0$ for all $q\in Y$. It is the restriction to $I$ of the 
left multiplication by another element $v\in Q$. Then $vu=e$ and $vY=0$ by the 
choice of $v$. The equalities $uv=e$ and $ev=ve=v$ follow again from the fact 
that $I$ has zero left annihilator. Thus $u,v\in eQe$, and $v$ is the inverse 
of $u$ in the ring $eQe$, as asserted in Claim 4.

\medbreak
Now we are able to continue as in the final part of the proof of \cite{Sk-Oy06, 
Lemma 6.4}. As explained there in detail, Claims 1 and 4 imply that each right 
ideal of $Q$ has the form $aQ+K$ where $a\in Q$ is an idempotent and $K$ is a 
nil right ideal of $Q$. Since the ring $Q$ is right Goldie, it has a largest 
nilpotent ideal $N$, and each nil right ideal of $Q$ is contained in $N$. In 
particular, $K\sbs N$. This shows that each right ideal of the factor ring 
$Q/N$ is generated by an idempotent. Hence $Q/N$ is semisimple Artinian, and 
$N$ is the Jacobson radical of $Q$.
\endproof

\section
3. Semiprimary Hopf module algebras

We wish to know that the ring $Q$ in Proposition 2.2 is actually a classical 
quotient ring of $A$. In spite of an effort made in \cite{Sk-Oy06} this 
conclusion in its full generality remains unproved. We will verify it under 
the additional assumption that the action of $H$ on $A$ is locally finite.

We will need only the properties of $Q$ established in Proposition 2.2, while 
the precise choice of the filter $\E'_H$ will not be significant any longer. 
Therefore our assumptions about $Q$ in this section are slightly more general. 
The final conclusion is presented in Proposition 3.9.

The first important step, where the local finiteness comes into play, consists 
in decomposing $Q$ as a direct product of $H$-simple algebras. This is done 
in Proposition 3.2, and here we have to apply one earlier result on the 
freeness of equivariant modules.

Let $B$ be any $H$-module algebra. A right $B$-module $M$ is said to be 
\emph{$H$-equivariant} if $M$ is equipped with a left $H$-module structure 
such that
$$
h(va)=\sum\,(h\1v)(h\2a)\quad\hbox{for all $h\in H$, $v\in M$, $a\in B$}.
$$
Suppose that $B$ is \emph{semilocal}, i.e., the factor ring of $B$ by the 
Jacobson radical is semisimple Artinian. Then the set $\,\Max B\,$ of maximal 
ideals of $B$ is finite, and the ring $B/P$ is simple Artinian for each 
$P\in\Max B$. If $M$ is a finitely generated right $B$-module, then $M/MP$ is 
a right $B/P$-module of finite length, and we define the \emph{rank} of $M$ at 
$P$ as
$$
r_P(M)={\lng M/MP\over\lng B/P}.
$$
In conformance with the notation of section 1 $P_H=\{a\in B\mid Ha\sbs P\}$ is 
the largest $H$-stable ideal of $B$ contained in $P$. Thus $P_H=0$ if and only 
if $P$ contains no nonzero $H$-stable ideals of $B$. Now let us recall 
\cite{Sk07, Lemma 7.5}:

\proclaim
Lemma 3.1.
Let $B$ be a semilocal $H$-module algebra and $M$ an $H$-equivariant finitely 
generated right $B$-module. Suppose that there exists $P\in\Max B$ such that 
$P_H=0$ and $r_P(M)\ge r_{P'}(M)$ for all $P'\in\Max B$. Then $M^n$ is a free 
$B$-module for some integer $n>0$.
\endproclaim

Lemma 3.1 can be applied, in particular, when $B$ is $H$-simple, but actually 
in some cases it can be used to show that $B$ is $H$-simple when this is not 
known beforehand. This is what we are going to do next.

\proclaim
Proposition 3.2.
Suppose that $Q$ is a semiprimary $H$-semiprime $H$-module algebra containing 
an $H$-stable subalgebra $A$ such that the action of $H$ on $A$ is locally 
finite and $I\cap A\ne0$ for each nonzero right ideal $I$ of $Q$. Then
there is an isomorphism of $H$-module algebras
$$
Q\cong Q_1\times\ldots\times Q_n
$$
where $Q_1,\ldots,Q_n$ are $H$-simple $H$-module algebras.
\endproclaim

\Proof.
Recall that every semiprimary ring satisfies DCC on finitely generated right 
ideals \cite{St, Ch. VIII, Prop. 5.5}. In particular, $Q$ has a minimal nonzero 
$H$-stable finitely generated right ideal, say $M$. For each $a\in A$ the 
$H$-stable right ideal $(Ha)Q$ is finitely generated by the local finiteness 
of the action on $A$. Since each nonzero right ideal of $Q$ contains nonzero 
elements of $A$, each nonzero $H$-stable right ideal of $Q$ contains therefore 
a nonzero $H$-stable finitely generated right ideal of $Q$.

It follows that $M$ is minimal in the set of all nonzero $H$-stable right 
ideals of $Q$. If $I$ is any two-sided $H$-stable ideal of $Q$, then either 
$MI=M$ or $MI=0$ since $MI$ is an $H$-stable right ideal and $MI\sbs M$.

Pick any maximal ideal $P\in\Max Q$ for which $r_P(M)$ attains the maximum 
value. Note that $M\ne MP$ since $M\ne0$. Then $MP_H\ne M$ too, and 
therefore $MP_H=0$. Thus $M$ is an $H$-equivariant finitely generated right 
$Q/P_H$-module and $P/P_H$ is a maximal ideal of the factor algebra $Q/P_H$ 
satisfying the hypotheses of Lemma 3.1. We conclude that a direct sum of 
several copies of $M$ is a free $Q/P_H$-module.

But then the assignment $I\mapsto MI$ gives an injection of the lattice of 
$H$-stable ideals of the $H$-module algebra $Q/P_H$ into the lattice of 
$H$-stable right ideals of $Q$ contained in $M$. Since $0$ and $M$ are the 
only two elements of the latter lattice, we deduce that the algebra $Q/P_H$ 
is $H$-simple.

Now $T=QM$ is a two-sided $H$-stable ideal of $Q$ such that $TP_H=0$. Since 
$(T\cap P_H)^2=0$, we must have $T\cap P_H=0$ by the $H$-semiprimeness of $Q$. 
On the other hand, $T+P_H=Q$ since this sum is an $H$-stable ideal of $Q$ 
properly containing $P_H$. The Chinese remainder theorem yields 
$Q\cong Q_1\times Q'$ where $Q_1=Q/P_H$ and $Q'=Q/T$.

Clearly $Q'$ is a semiprimary $H$-module algebra, and the projection $\pi$ of 
$Q$ onto $Q'$ is a homomorphism of $H$-module algebras. Hence $\pi(A)$ is 
an $H$-stable subalgebra of $Q'$ on which the action of $H$ is locally finite. 
Each $H$-stable right ideal of $Q'$ can be written as $\pi(J)$ where $J$ is an 
$H$-stable right ideal of $Q$ lying in the kernel of the other projection 
$Q\to Q_1$. If $J\ne0$, then $J\cap A\ne0$ by the hypothesis, and it follows 
that $\pi(J)\cap\pi(A)\ne0$ since the map $\pi|_J$ is injective. Also $Q'$ is 
$H$-semiprime since so is $Q$.

Thus $Q'$ satisfies the same assumptions as $Q$, but has fewer maximal 
ideals. We have seen that $Q_1$ is $H$-simple. Proceeding by induction on the 
cardinality of the set $\Max Q$, we may assume that $Q'$ is a direct product 
of finitely many $H$-simple $H$-module algebras, and the proof is completed.  
\endproof

\proclaim
Lemma 3.3.
The $H$-module algebra $Q$ in Proposition\/ {\rm3.2} is in fact 
$S(H)$-semiprime and each direct factor $Q_i$ is $S(H)$-simple.
\endproclaim

\Proof.
By an argument given in the proof of Proposition 3.2 any $H$-module algebra 
isomorphic to a direct factor of $Q$ satisfies the assumptions imposed on $Q$ 
in the statement of Proposition 3.2. Therefore it suffices to consider the 
case when $n=1$ and $Q$ is $H$-simple.

For any $S(H)$-stable ideal $I$ of $Q$ its left annihilator is an $H$-stable 
ideal of $Q$. Since $Q$ is $H$-simple, we must have $\lann_QI=0$ whenever 
$I\ne0$.  In particular, a nonzero $S(H)$-stable ideal cannot be nilpotent. So 
$Q$ is $S(H)$-semiprime.

Now we can apply Proposition 3.2, replacing $H$ with its Hopf subalgebra 
$S(H)$. It shows that $Q$ is isomorphic as an $S(H)$-module algebra to a 
direct product of several $S(H)$-simple $S(H)$-module algebras. The direct 
factors may be identified with minimal nonzero $S(H)$-stable ideals of $Q$. If 
$I_1$, $I_2$ are two different such ideals, then $I_1I_2\sbs I_1\cap I_2=0$, 
but this is impossible since $\,\lann_QI_2=0$, as we have seen already. Hence 
$Q$ is $S(H)$-simple.
\endproof

Let $B$ be an arbitrary $H$-module algebra. For a right $H$-comodule $U$ and a 
right $B$-module $V$ we will consider the vector space $U\ot V$ as a right 
$B$-module with respect to the \emph{twisted action} of $B$ defined by the 
rule
$$
(u\ot v)a=\sum u\0\ot v\bigl(S(u\1)a\bigr),\qquad u\in U,\ v\in V,\ a\in B,
$$
where $\sum u\0\ot u\1\in U\ot H$ is the symbolic notation for the image of 
$u$ under the comodule structure map $U\to U\ot H$.

We will also need similar tensoring operations on the left modules. Given a 
left $B$-module $V$ and $U$ as above, there is a left $B$-module structure 
on the vector space $V\ot U$ defined by the rule
$$
a(v\ot u)=\sum(u\1a)v\ot u\0,\qquad u\in U,\ v\in V,\ a\in B.
$$

\proclaim
Lemma 3.4.
Suppose that $\G$ is a right Gabriel topology on $A$ such that all elements 
of $S(H)$ operate on $A$ as $\G$-continuous transformations. If $V$ is a 
$\G$-torsion right $A$-module{\rm,} then so is $U\ot V$ for any right 
$H$-comodule $U$.
\endproclaim

\Proof.
Let $u\in U$ and $v\in V$. We have $vI=0$ for some $I\in\G$. Since 
$\sum u\0\ot u\1\in U\ot C$ for some $C\in\F$, it follows from the formula for 
the action of $A$ in $U\ot V$ that $u\ot v$ is annihilated by the right ideal 
$I_{S(C)}$ of $A$. But $I_{S(C)}\in\G$, and therefore $u\ot v$ lies in the 
$\G$-torsion submodule of $U\ot V$. Since such elements span the whole 
$U\ot V$, the conclusion follows.
\endproof

\proclaim
Lemma 3.5.
Suppose that $B$ is an $S(H)$-simple $H$-module algebra. If $K$ is a simple 
right ideal of $B$ and $V$ is any nonzero right $B$-module{\rm,} then $K$ 
embeds in the right $B$-module $U\ot V$ for some finite dimensional right 
$H$-comodule $U$.
\endproclaim

\Proof.
We may regard $H$ as a right $H$-comodule with respect to the comultiplication 
in $H$. If $a\in B$ annihilates $H\ot V$, then
$$
\sum h\1\ot v\bigl(S(h\2)a\bigr)=0\quad
\hbox{for all $h\in H$ and $v\in V$},
$$
and applying the map $\ep\ot\id:H\ot V\to V$, we get $v\bigl(S(h)a\bigr)=0$. 
Hence $S(H)a$ is contained in the annihilator $I$ of the $B$-module $V$, i.e., 
$a$ lies in the largest $S(H)$-stable ideal $I_{S(H)}$ of $B$ contained in $I$. 
Since $V\ne0$, we have $I\ne B$, but then $I_{S(H)}=0$ by the $S(H)$-simplicity 
of $B$.

This shows that $H\ot V$ is a faithful $B$-module. Therefore there exists an 
element $t\in H\ot V$ such that $tK\ne0$. We have $tK\cong K$ since $K$ is a 
simple right $B$-module, and $t\in C\ot V$ for some $C\in\F$ since $H$ is the 
union of finite dimensional subcoalgebras. Then $tK\sbs C\ot V$, and we 
may take $U=C$, a subcomodule of $H$.
\endproof

\proclaim
Lemma 3.6.
Suppose that $Q$ is a semiprimary $H$-semiprime $H$-module algebra containing 
an $H$-stable subalgebra $A$ on which the action of $H$ is locally finite. 
Suppose also that $\G$ is a right Gabriel topology on $A$ such that all elements 
of $S(H)$ operate on $A$ as $\G$-continuous transformations and the following 
two properties hold\/{\rm:}

\item(a)
$\lann_QI=0$ for each $I\in\G,$

\item(b)
for each $q\in Q$ there exists $I\in\G$ such that $qI\sbs A$.

Then all right $Q$-modules are $\G$-torsionfree as right $A$-modules{\rm,} and 
therefore $IQ=Q$ for each $I\in\G$.
\endproclaim

\Proof.
If in (b) $q\ne0$, then $qI\ne0$ by (a). This shows that each nonzero right 
ideal of $Q$ has nonzero intersection with $A$. So the assumptions of 
Proposition 3.2 are satisfied, and we conclude that 
$Q\cong Q_1\times\ldots\times Q_n$ where $Q_1,\ldots,Q_n$ are $S(H)$-simple 
$H$-module algebras by Lemma 3.3.

By (a) $Q$ is a $\G$-torsionfree right $A$-module. Hence so are all right 
ideals of $Q$. Let $M$ be any right $Q$-module. We have 
$M\cong M_1\oplus\ldots\oplus M_n$ where $M_i$, for each $i$, is a $Q_i$-module 
on which $Q$ acts via the projection $Q\to Q_i$. To prove that $M$ is 
$\G$-torsionfree it suffices to consider the case when $M=M_i$ for some $i$.

Denote by $N$ the largest $\G$-torsion $A$-submodule of $M$. If $x\in N$ and 
$q\in Q$, then the coset $xq+N$ is annihilated by a right ideal in $\G$, 
according to (b). Since the $A$-module $M/N$ is $\G$-torsionfree, it follows 
that $Nq\sbs N$ for each $q\in Q$. In other words, $N$ is a $Q$-submodule of 
$M$. Assuming $M$ to be a $Q_i$-module, we conclude that so is $N$.

Since $Q_i$ is semiprimary, it has a simple right ideal, say $K$. If $N\ne0$, 
then, by Lemma 3.5, $K$ embeds in $U\ot N$ for some right $H$-comodule $U$. In 
this case $K$ has to be $\G$-torsion by Lemma 3.4. But this is impossible 
since $K$ is isomorphic to a right ideal of $Q$, and therefore $K$ is 
$\G$-torsionfree, as we have observed already. Thus $N=0$, and $M$ is indeed 
$\G$-torsionfree.

In particular, the right $Q$-module $Q/IQ$ is $\G$-torsionfree for any right 
ideal $I$ of $A$. On the other hand, $Q/IQ$ is $\G$-torsion whenever $I\in\G$ 
since the right $A$-modules $Q/A$ and $A/I$ are $\G$-torsion. In this case we 
must have $\,Q/IQ=0$.
\endproof

The conclusion of Lemma 3.6 implies that $Q$ is a perfect right localization 
of $A$ (see \cite{St, Ch. XI, Th. 2.1}). However, the final goal has not been 
reached yet.

By Lemma 3.6 no nonzero element of a right $Q$-module is annihilated by a 
right ideal in $\G$. We will need a similar conclusion for left $Q$-modules. 
This will require more delicate arguments since $\G$ is not a left Gabriel 
topology. For a left $A$-module $M$ put
$$
T_\G(M)=\{x\in M\mid x\hbox{ is annihilated by a right ideal in $\G$}\}.
$$
For each right ideal $I\in\G$ the set $\Ann_MI=\{x\in M\mid Ix=0\}$ is a 
submodule of $M$. Since $\G$ is a filter, the set of all such submodules is 
directed by inclusion. Hence
$$
T_\G(M)=\bigcup_{I\in\G}\Ann_MI
$$
is a submodule too. However, we cannot be sure that $T_\G(M)$ is a 
$Q$-submodule when $M$ is a left $Q$-module.

Note that $T_\G(M)$ is stable under all endomorphisms of $M$. In particular, 
if $D$ is a skew field contained in the endomorphism ring 
$\,\End_A\mskip-2mu M$, then $T_\G(M)$ is a vector space over $D$.

\proclaim
Lemma 3.7.
In addition to the hypothesis of Lemma\/ {\rm3.6} assume that $\,\rann_AI=0$ 
for each $I\in\G$. Then $T_\G(M)=0$ for each left $Q$-module $M$.
\endproclaim

\Proof.
If $I\in\G$, then $A\cap\rann_QI=\rann_AI=0$. Since $\rann_QI$ is a right 
ideal of $Q$ having zero intersection with $A$, we must have $\rann_QI=0$. 
This shows that $T_\G(Q)=0$, and therefore $T_\G(L)=0$ for each left ideal $L$ 
of $Q$.

It is clear that $T_\G$ is a left exact functor. Thus for each exact sequence 
of left $Q$-modules $0\to M'\to M\to M''\to0$ there is an exact sequence 
of left $A$-modules
$$
0\to T_\G(M')\to T_\G(M)\to T_\G(M''),\eqno(*)
$$
and it follows that $T_\G(M)=0$ whenever $T_\G(M')=0$ and $T_\G(M'')=0$. 
Since the Jacobson radical $J$ of $Q$ is nilpotent, each left $Q$-module has a 
finite chain of submodules with factors annihilated by $J$. Hence it suffices 
to prove that $T_\G(M)=0$ when $JM=0$. Since the factor ring $Q/J$ is 
semisimple Artinian, any such a module $M$ is semisimple. Since $T_\G$ is 
an additive functor, the proof of the equality $T_\G(M)=0$ reduces further to 
the case when $M$ is simple.

There are finitely many isomorphism classes of simple $Q$-modules. Let 
$V_1,\ldots,V_p$ be a full set of pairwise nonisomorphic simple left 
$Q$-modules. For each $i$ the endomorphism ring $E_i=\End_QV_i$ is a skew 
field and $V_i$ is a finite dimensional vector space over $E_i$. Then 
$T_\G(V_i)$ is a vector subspace of $V_i$. Put
$$
\mu=\max_{i=1,\ldots,p}{\dim_{E_i}\!T_\G(V_i)\over\dim_{E_i}\!V_i}.
$$
We will show that $\mu=0$. This will yield $T_\G(V_i)=0$ for all $i$, and the 
proof of the lemma will be completed.

\proclaim
Claim 1.
Suppose that $M$ is a left $Q$-module and $D$ is a skew field contained in the 
endomorphism ring $\,\End_QM$. If\/ $\dim_DM<\infty,$ then
$$
\dim_DT_\G(M)\le\mu\cdot\dim_DM.
$$
\endproclaim

Consider first the case when $M$ is an isotypic semisimple left $Q$-module. 
In other words, $M$ is a direct sum of a possibly infinite family of copies of 
some simple module $V_i$. With $S=\Hom_Q(V_i,M)$ we have 
$M\cong S\ot_{E_i}\!V_i$ as $Q$-modules and as $\End_QM$-modules. Here $Q$ 
acts on $V_i$, while $\End_QM$ acts on $S$. It follows that
$$
T_\G(M)\cong S\ot_{E_i}\!T_\G(V_i)
$$
as $\,\End_QM$-modules. We get $\,\dim_DM=(\dim_DS)(\dim_{E_i}\!V_i)\,$ and
$$
\dim_DT_\G(M)=(\dim_DS)\bigl(\dim_{E_i}\!T_\G(V_i)\bigr).
$$
The assumption $\dim_DM<\infty$ implies that $\dim_DS<\infty$, and the claim 
follows from the inequality
$$
\dim_{E_i}\!T_\G(V_i)\le\mu\cdot\dim_{E_i}\!V_i\,.
$$

In the general case we proceed as follows. If $M'$ is any submodule of $M$ 
stable under the action of $D$, then $D$ embeds in $\,\End_QM'$ and in 
$\,\End_QM''$ where we put $M''=M/M'$. Then
$$
\openup1\jot
\eqalign{
\dim_DM&{}=\dim_DM'+\dim_DM'',\cr
\dim_DT_\G(M)&{}\le\dim_DT_\G(M')+\dim_DT_\G(M'')
}\eqno(**)
$$
where the last inequality follows from the exact sequence of $D$-vector spaces 
$(*)$. If Claim 1 is true for the $Q$-modules $M'$ and $M''$, it is clear that 
Claim 1 is true for $M$ as well.

We can use this argument with $M'=JM$. Since $J$ is nilpotent, verification of 
Claim 1 is thus reduced to the case when $JM=0$, and so $M$ is semisimple. But 
then the isotypic components of $M$ are stable under all endomorphisms, and 
the proof reduces in a similar way to the case considered at the beginning.  
Thus Claim 1 has been verified.

\proclaim
Claim 2.
If in Claim\/ {\rm1} the equality $\,\dim_DT_\G(M)=\mu\cdot\dim_DM$ is 
attained{\rm,} then
$$
\dim_DT_\G(M')=\mu\cdot\dim_DM'
$$
for any $Q$-submodule $M'$ of $M$ stable under the action of $D$.
\endproclaim

Indeed, $\,\dim_DT_\G(M')\le\mu\cdot\dim_DM'$ and 
$\,\dim_DT_\G(M'')\le\mu\cdot\dim_DM''$ where $M''=M/M'$. If one of these two 
inequalities were strict, then we would get
$$
\dim_DT_\G(M)<\mu\cdot\dim_DM
$$
from $(**)$, a contradiction. Thus both inequalities are in fact equalities.

\proclaim
Claim 3.
Suppose that $B$ is an $S(H)$-simple $H$-module algebra{\rm,} $L$ is a simple 
left ideal of $B,$ and $V$ is any nonzero left $B$-module. There exists 
a finite dimensional right $S(H)$-comodule $U$ such that $L$ embeds as a 
$B$-submodule in $V\ot U$ with the twisted left $B$-module structure.  
\endproclaim

This is an analog of Lemma 3.5 with a similar proof. Considering $S(H)$ as a 
right $H$-comodule with respect to the comultiplication in $H$, the twisted 
left $B$-module $V\ot S(H)$ is faithful. Hence $L$ embeds in $V\ot S(H)$, and 
therefore in $V\ot U$ where $U=S(C)$ for some finite dimensional subcoalgebra 
$C$ of $H$.

\medbreak
We are now in a position to complete the proof of Lemma 3.7. Among the simple 
left $Q$-modules $V_1,\ldots,V_p$ pick $V_j$ such that
$$
\dim_{E_j}\!T_\G(V_j)=\mu\cdot\dim_{E_j}\!V_j.
$$
Let $Q_1,\ldots,Q_n$ be the $H$-simple direct factors of $Q$ given by 
Proposition 3.2. Then $V_j$ is a $Q_i$-module for some $i$. Since $Q_i$ is 
semiprimary, it has a simple left ideal, say $L$. By Lemma 3.3 $Q_i$ is 
$S(H)$-simple. Therefore Claim 3 shows that $L$ is isomorphic to a submodule 
of the twisted $Q_i$-module $M=V_j\ot U$ for some finite dimensional right 
$S(H)$-comodule $U$.

Now note that $T_\G(V_j)\ot U\sbs T_\G(M)$. Indeed, if $v\in T_\G(V_j)$ and 
$u\in U$, then $Iv=0$ for some $I\in\G$ and $\sum u\0\ot u\1\in U\ot S(C)$ 
for some $C\in\F$. It follows then from the formula for the twisted action of 
$Q$ in $V_j\ot U$ that $v\ot u$ is annihilated by the right ideal $I_{S(C)}$ 
of $A$. Since $I_{S(C)}\in\G$, we get $v\ot u\in T_\G(M)$.

The skew field $E_j=\End_QV_j$ embeds in $\End_QM$ in a natural way, and
$$
\openup1\jot
\eqalign{
\dim_{E_j}\!M&{}=(\dim_{E_j}\!V_j)(\dim_kU),\cr
\dim_{E_j}\!T_\G(M)&{}\ge\bigl(\dim_{E_j}\!T_\G(V_j)\bigr)(\dim_kU).
}
$$
Hence $\,\dim_{E_j}\!T_\G(M)\ge\mu\cdot\dim_{E_j}\!M$. By Claim 1 the opposite 
inequality is also true, and so we must have an equality here. But then Claim 2 
shows that
$$
\dim_{E_j}\!T_\G(M')=\mu\cdot\dim_{E_j}\!M'
$$
for each $Q$-submodule $M'$ of $M$ stable under the action of $E_j$.

Take $M'$ to be the sum of all $Q$-submodules of $M$ isomorphic to $L$. 
Obviously, $M'$ is stable under all endomorphisms of $M$, and so the previous 
equality must hold. Note that $M'\ne0$. On the other hand, $L$ is isomorphic 
to a right ideal of $Q$. As we have seen, this entails $T_\G(L)=0$. Since $M'$ 
is an isotypic semisimple $Q$-module, we get $T_\G(M')=0$. It follows that 
$\mu=0$, and we are done.
\endproof

\proclaim
Lemma 3.8.
Let $R$ be a semiprime right Goldie subring of a semisimple Artinian ring $S$. 
Suppose that $\G$ is a set of right ideals of $R$ with the two properties\/{\rm:}

\item(a)
$\lann_SI=0$ for each $I\in\G,$

\item(b)
for each $x\in S$ there exists $I\in\G$ such that $xI\sbs R$.

Then $S$ is a classical right quotient ring of $R,$ and each right ideal 
$I\in\G$ contains a regular element of $R$.
\endproclaim

\Proof.
By the Goldie theorem the ring $R$ has a semisimple Artinian classical right 
quotient ring $Q$. It is known from the proof of Goldie's theorem that $Q$ is 
the localization of $R$ with respect to the filter $\E$ of all essential right 
ideals of $R$, and a right ideal of $R$ is essential if and only if it 
contains a regular element of $R$.

Considering $S$ as a right $R$-module, denote by $T$ the set of all elements 
$x\in S$ whose right annihilator in $R$ belongs to $\E$. Then $T$ is a right 
$R$-submodule of $S$, and $T\cap R=0$ since no nonzero element of $R$ is 
annihilated by a regular element of $R$. On the other hand, (a) and (b) imply 
that each nonzero right $R$-submodule of $S$ has nonzero intersection with 
$R$. Hence $T=0$.

Thus we have shown that $\lann_SI=0$ for each $I\in\E$. In other words, 
$\lann_Su=0$ for each regular element $u$ of $R$, so that all regular elements 
of $R$ are left regular in $S$. But left regular elements of a right or left 
Artinian ring are invertible. We conclude that regular elements of $R$ are 
invertible in $S$.

By the universality property of the Ore localizations, the embedding $R\to S$ 
extends to a ring homomorphism $\ph:Q\to S$, and $\ph$ is injective since 
$\Ker\ph\cap R=0$. So $Q$ is identified with a subring of $S$.

Let $I\in\G$. By (a) $I$ has zero left annihilator in $Q$. Since $Q$ is 
semisimple Artinian, its right ideal $IQ$ is generated by an idempotent, say 
$e$. Noting that $(1-e)I=0$, we deduce that $e=1$. Thus $IQ=Q$, which means 
precisely that $I$ contains a regular element of $R$.

It follows now from (b) that for each $x\in S$ we have $xu\in R$ for some 
regular element $u$ of $R$, whence $x=(xu)u^{-1}\in Q$. We conclude that $S=Q$.
\endproof

\proclaim
Proposition 3.9.
Suppose that $Q$ is a semiprimary $H$-semiprime $H$-module algebra containing 
a right Noetherian $H$-stable subalgebra $A$ on which the action of $H$ is 
locally finite. Suppose also that $\G$ is a right Gabriel topology on $A$ 
such that all elements of $S(H)$ operate on $A$ as $\G$-continuous 
transformations and the following two properties hold\/{\rm:}

\item(a)
$\lann_QI=0$ and\/ $\rann_AI=0$ for each $I\in\G,$

\item(b)
for each $q\in Q$ there exists $I\in\G$ such that $qI\sbs A$.

Then $Q$ is a classical right quotient ring of $A$.
\endproclaim

\Proof.
We proceed in several steps.

\proclaim
Claim 1.
If $M$ is any maximal ideal of $Q,$ then $M\cap A$ is a prime ideal of $A$.
\endproclaim

Recall that for any nonzero right module over a right Noetherian ring the set 
of annihilators of nonzero submodules contains a prime ideal. For example, all 
maximal elements of this set are prime ideals \cite{Goo-W, Prop. 3.12}. Let 
$P$ be any prime annihilator of a nonzero submodule of the right $A$-module 
$Q/M$. Thus $P$ is a prime ideal of $A$ such that $M\cap A\sbs P$ and 
$yP\sbs M$ for some $y\in Q$, $y\notin M$.

For each $I\in\G$ put
$$
L(I)=\{x\in Q\mid xIP\sbs M\}.
$$
It is clear that $L(I)$ is a left ideal of $Q$ and $M\sbs L(I)$, so that 
$L(I)/M$ is a left ideal of the factor ring $Q/M$. Since the simple ring $Q/M$ 
is Artinian, the set of left ideals $\{L(I)\mid I\in\G\}$ has a maximal element, 
say $L_0$. Moreover, since the correspondence $I\mapsto L(I)$ reverses 
inclusions, this set is directed by inclusion, whence $L_0$ is in fact its 
largest element. Thus $L(I)\sbs L_0$ for each $I\in\G$.

Note that $L(A)=\{x\in Q\mid xP\sbs M\}$ is a right $A$-submodule of $Q$. Let 
$x\in L(A)$ and $q\in Q$. By (b) there exists $I\in\G$ such that $qI\sbs A$. 
Then
$$
xqIP\sbs xAP\sbs xP\sbs M,
$$
yielding $xq\in L(I)\sbs L_0$. This shows that $L(A)Q\sbs L_0$. Note that 
$L(A)Q$ is a two-sided ideal of $Q$ containing $M$. But $L(A)\not\sbs M$ by 
the choice of $P$, whence $L(A)Q\not\sbs M$. Since $Q/M$ is a simple ring, we 
must have $L(A)Q=Q$. It follows that $L_0=Q$ too.

Now pick $I\in\G$ such that $L_0=L(I)$. Since $1\in L_0$, we get $IP\sbs M$. 
This means that the image of $P$ in $Q/M$ is contained in the left 
$A$-submodule $T_\G(Q/M)$, in the notation of Lemma 3.7. By that lemma 
$T_\G(Q/M)=0$, which entails $P\sbs M$. Hence $P=M\cap A$, and Claim 1 is thus 
proved.

\proclaim
Claim 2.
Denote by $N$ the prime radical of $A$ and by $J$ the Jacobson radical of $Q$. 
Then $N=J\cap A$. 
\endproclaim

Let $M_1,\ldots,M_k$ be all the maximal ideals of $Q$. Then $J=\bigcap M_i$. 
By Claim 1 $P_i=M_i\cap A$ is a prime ideal of $A$ for each $i$. We have 
$J\cap A=N'$ where $N'=\bigcap P_i$. Since $J$ is nilpotent, $N'$ is a 
nilpotent ideal of $A$. Hence $N'$ is contained in each prime ideal of $A$, 
and therefore $N'=N$.

\proclaim
Claim 3.
The factor ring $Q/J$ is a classical right quotient ring of $A/N$.
\endproclaim

The ring $S=Q/J$ is semisimple Artinian, while $R=A/N$ is semiprime right 
Noetherian. Let $\pi:Q\to S$ be the canonical homomorphism. By Claim 2 $R$ is 
identified with the subring $\pi(A)$ of $S$.

Consider the set $\{\pi(I)\mid I\in\G\}$ of right ideals of $R$.
If $q\in Q$, then it follows from (b) that $\pi(q)\pi(I)\sbs\pi(A)$ for some 
$I\in\G$. For each $I\in\G$ we have $IQ=Q$ by Lemma 3.6, whence $\pi(I)S=S$, 
and therefore $\pi(I)$ has zero left annihilator in $S$. Thus we meet the 
hypothesis of Lemma 3.8, and Claim 3 follows.

\medbreak
Denote by $\C$ the set of all elements $u\in A$ which are regular modulo $N$, 
i.e., whose images $\pi(u)$ in the ring $\pi(A)\cong A/N$ are regular elements 
of that ring. If $u\in\C$, then $\pi(u)$ is invertible in $Q/J$ by Claim 3; 
this implies that $u$ is invertible in $Q$ since $J$ is the Jacobson radical 
of $Q$. This shows that all elements of $\C$ are regular in $A$. Conversely, 
each regular element of a right Noetherian ring is regular modulo the prime 
radical (see \cite{Goo-W, Lemma 11.8} or \cite{Mc-R, 4.1.3}). Thus $\C$ is the 
set of all regular elements of $A$.

By Lemma 3.8 each right ideal in the set $\{\pi(I)\mid I\in\G\}$ contains a 
regular element of the ring $\pi(A)$. Therefore $I\cap\C\ne\varnothing$ for 
all $I\in\G$. It follows now from condition (b) that for each $q\in Q$ there 
exists $u\in\C$ such that $qu\in A$. We have seen already that all elements of 
$\C$ are invertible in $Q$. These properties characterize $Q$ as a classical 
right quotient ring of $A$.
\endproof

\section
4. Final results

In the first result of this section we complete the work on Noetherian 
$H$-semiprime $H$-module algebras done in the preceding sections. This will 
then be used to derive results on bijectivity of the antipode and on flatness 
over coideal subalgebras.

\proclaim
Theorem 4.1.
Let $A$ be a right Noetherian $H$-semiprime $H$-module algebra such that the 
action of $H$ on $A$ is locally finite. Then $A$ has a right Artinian 
classical right quotient ring.
\endproclaim

\Proof.
Let $Q$ be the quotient ring of $A$ with respect to the filter of right 
ideals $\G=\E'_H$. By Proposition 1.4 $\G$ is a Gabriel topology and $A$ is 
$\G$-torsionfree. Since $Q$ is an essential extension of $A$ in the category 
of right $A$-modules, $Q$ is $\G$-torsionfree as well. Combined with Corollary 
1.10 this amounts to condition (a) in the statement of Proposition 3.9. 
Condition (b) is satisfied by the construction of $Q$. By Proposition 2.2 $Q$ 
is semiprimary and $H$-semiprime. Thus the assumptions of Proposition 3.9 are 
fulfilled, and therefore $Q$ is a classical right quotient ring of $A$. Since 
$A$ is right Noetherian, so too is $Q$. Since $Q$ is also semiprimary, it has 
to be right Artinian (see \cite{St, Ch. VIII, Prop. 1.12}).
\endproof

From \cite{Sk11, Th. 1.1} we deduce that the quotient ring $Q$ in Theorem 4.1 
is quasi-Frobenius, but this fact will not be needed.

\smallskip
The dual Hopf algebra $\Hd$ consists of all linear functions $H\to k$ 
vanishing on an ideal of finite codimension in $H$. There is an action of 
$\Hd$ on $H$ defined by the rule
$$
f\rhu h=\sum\, f(h\2)\,h\1\,,\qquad f\in\Hd,\ \,h\in H.
$$
It makes $H$ into a left $\Hd$-module algebra. Right coideals of $H$ are 
stable under this action of $\Hd$. Since each element of $H$ is contained in a 
finite dimensional subcoalgebra, the action of $\Hd$ on $H$ is locally finite.

\proclaim
Lemma 4.2.
Suppose that $H$ is a residually finite dimensional Hopf algebra. Then $H$ is 
an $\Hd$-simple $\Hd$-module algebra. Each right coideal subalgebra $A$ of $H$ 
is an $\Hd$-prime $\Hd$-module algebra.
\endproclaim

\Proof.
Since $H$ is residually finite dimensional, the $\Hd$-submodules of $H$ are 
precisely the right coideals. Now, if $I$ is a right ideal of $H$ such that 
$\De(I)\sbs I\ot H$, then $I$ may be regarded as a Hopf module. By the 
structure of Hopf modules (see \cite{Mo, 1.9.4})
$$
I=I^{\co H}H\quad{\rm where}\ \,I^{\co H}=\{h\in I\mid\De(h)=h\ot1\}.
$$
Since $H^{\co H}=k$, we deduce that $I^{\co H}$ equals either 0 or $k$. Hence 
either $I=0$ or $I=H$. In other words, $0$ and $H$ are the only two 
$\Hd$-stable right ideals of $H$. In particular, $H$ is an $\Hd$-simple.

Any right coideal subalgebra $A$ of $H$ is an $\Hd$-stable subalgebra, and so 
is itself an $\Hd$-module algebra. Suppose that $I$ is a nonzero $\Hd$-stable 
ideal of $A$. Then $IH$ is a nonzero $\Hd$-stable right ideal of $H$, whence 
$IH=H$. If $J$ is another nonzero $\Hd$-stable ideal of $A$, then $JH=H$ too. 
It follows that $IJH=H$, and therefore $IJ\ne0$. Thus $A$ is $\Hd$-prime.
\endproof

\proclaim
Theorem 4.3.
Let $H$ be either right or left Noetherian residually finite dimensional Hopf 
algebra. Then its antipode $S:H\to H$ is bijective. Hence $H$ is right and 
left Noetherian simultaneously.
\endproclaim

\Proof.
According to \cite{Sk06, Th. A} $S$ is bijective whenever $H$ can be embedded 
into a left perfect ring $Q$ such that $Q$ is an essential extension of $H$ as 
a right $H$-module. By Lemma 4.2 $H$ is $\Hd$-simple, and therefore 
$\Hd$-semiprime, as an $\Hd$-module algebra. If $H$ is right Noetherian, the 
required embedding is provided already by Proposition 2.2 (semiprimary rings 
are left perfect).

If $H$ is left Noetherian, we consider the Hopf algebra $H\opcop$ obtained 
from $H$ by taking the opposite multiplication and comultiplication. It has 
the same antipode $S$, but is right Noetherian. So we can refer to the case 
already treated.

Bijectivity of $S$ implies that $S$ is an antiautomorphism of $H$ as a Hopf 
algebra. In particular, $H\cong H\op$ as algebras. Therefore the right hand 
properties of $H$ are equivalent to the left hand ones.
\endproof

In the next theorem we repeat results which can be found in \cite{Sk10, Th. 
1.8, Cor. 1.9}. However, we present only the part concerned with flatness and 
provide a proof which bypasses the category equivalences considered in 
\cite{Sk10}.

\proclaim
Theorem 4.4.
Let $A$ be a right coideal subalgebra of a residually finite dimensional Hopf 
algebra $H$. Suppose that $A$ and $H$ have right Artinian classical right 
quotient rings $Q(A)$ and $Q(H)$. Then $H$ is left $A$-flat. Moreover{\rm,} if 
$A$ is a Hopf subalgebra{\rm,} then $H$ is left and right faithfully $A$-flat.  
\endproclaim

\Proof.
First we note that the embedding $H\hrar Q(H)$ enables us to apply \cite{Sk06, 
Th. A} and conclude that the antipode $S:H\to H$ is bijective. This fact will 
be used without further notice. Next, the action of $\Hd$ on $A$ and $H$ 
extends to $Q(A)$ and $Q(H)$ by \cite{Sk-Oy06, Th. 2.2}.

\proclaim
Claim 1.
$Q(A)$ is an $\Hd$-simple $\Hd$-module algebra.
\endproclaim

If $I$ and $J$ are two nonzero $\Hd$-stable ideals of $Q(A)$, then $I\cap A$ 
and $J\cap A$ are nonzero $\Hd$-stable ideals of $A$. Since $A$ is $\Hd$-prime 
by Lemma 4.2, we deduce that $IJ\ne0$. Thus $Q(A)$ is $\Hd$-prime.

Since the action of $\Hd$ on $A$ is locally finite, we can now apply 
Proposition 3.2. It shows that $Q(A)$ is a direct product of finitely many 
$\Hd$-simple $\Hd$-module algebras $Q_1,\ldots,Q_n$. Since $Q(A)$ is 
$\Hd$-prime, we cannot have $n>1$. Hence $n=1$, and so $Q(A)=Q_1$ is indeed 
$\Hd$-simple.

\proclaim
Claim 2.
The inclusion $A\hrar H$ extends to a ring homomorphism $Q(A)\to Q(H)$.
\endproclaim

This is a special case of \cite{Sk10, Lemma 1.7}. We do not offer any 
improvements in its proof.

\medbreak
Let further $V$ be a right $A$-module, and put $M=V\ot H$. We will view $M$ as 
an $\Hd$-equivariant right $A$-module with the actions of $A$ and $\Hd$ defined 
as follows:
$$
(v\ot h)\,a=\sum va\1\ot ha\2,\,\qquad f\rhu(v\ot h)=v\ot(f\rhu h)
$$
for $v\in V$, $h\in H$, $a\in A$ and $f\in\Hd$.

In this way $?\ot H$ becomes 
a functor from the category of right $A$-modules to the category of 
$\Hd$-equivariant right $A$-modules. Note that the action of $\Hd$ on $M$ is 
locally finite since so is the action of $\Hd$ on $H$.

\proclaim
Claim 3.
The right $Q(A)$-module $M\ot_AQ(A)$ is projective. Moreover{\rm,} 
$M^n\ot_AQ(A)$ is a free $Q(A)$-module for some integer $n>0$.
\endproclaim

Since $Q(A)$ is an extension of $A$ in the category of $\Hd$-module algebras, 
there is a well-defined action of $\Hd$ on $M\ot_AQ(A)$ such that
$$
f\rhu(x\ot q)=\sum\,(f\1\rhu x)\ot(f\2\rhu q)
\quad\hbox{for $f\in\Hd$, $x\in M$, $q\in Q(A)$}.
$$
It makes $M\ot_AQ(A)$ an $\Hd$-equivariant right $Q(A)$-module. If $U$ is any 
$\Hd$-stable subspace of $M$, then the $A$-submodule $UA$ generated by $U$ is 
$\Hd$-stable too, and
$$
F_U=(UA)\ot_AQ(A)
$$
is an $\Hd$-equivariant right $Q(A)$-module which may be identified with a 
submodule of $M\ot_AQ(A)$ since $Q(A)$ is left $A$-flat by the standard 
properties of classical quotient rings. If $\dim U<\infty$, then $F_U$ is 
finitely generated. Each $\Hd$-equivariant finitely generated right 
$Q(A)$-module is projective since $Q(A)$ is $\Hd$-simple, so that Lemma 3.1 can 
be applied with a suitable choice of the maximal ideal $P$. Thus
$$
\{F_U\mid\hbox{$U$ is a finite dimensional $\Hd$-submodule of $M$}\}
$$
is a directed set of projective submodules of the right $Q(A)$-module 
$M\ot_AQ(A)$. The union of this family of submodules gives the whole module 
since the action of $\Hd$ on $M$ is locally finite. We arrive at the first 
conclusion of Claim 3, observing that inductive direct limits of flat modules 
are flat and that all flat right modules over a right Artinian ring are 
projective.

Take $n$ to be the greatest common divisor of the lengths of simple factor 
rings of the right Artinian ring $Q(A)$. It follows from Lemma 3.1 and the 
Krull-Schmidt Theorem that for each $\Hd$-equivariant finitely generated right 
$Q(A)$-module $K$ the $Q(A)$-module $K^n$ is free for exactly this value of 
$n$ which does not depend on $K$. A basis of $M^n\ot_AQ(A)$ over $Q(A)$ can 
then be obtained by a suitable application of Zorn's Lemma (see \cite{Sk07, 
Th. 7.6}).

\proclaim
Claim 4.
There is an isomorphism of right $H$-modules $M\ot_AH\cong(V\ot_AH)\ot H$ 
where $H$ is assumed to act by right multiplications on the last tensorands.
\endproclaim

We have $M\ot_AH\cong(V\ot H\ot H)/R$ where $R$ is the subspace of 
$V\ot H\ot H$ spanned by
$$
\{\sum va\1\ot ga\2\ot h-v\ot g\ot ah\mid v\in V,\ a\in A,\ g,h\in H\}.
$$
Denote by $R'$ the subspace of $V\ot H\ot H$ spanned by
$$
\{va\ot g\ot h-\sum v\ot gS^{-1}(a\2)\ot a\1h\mid v\in V,\ a\in A,\ g,h\in H\}.
$$
In fact $R'\sbs R$ since
$$
va\ot g\ot h=\sum va\1\ot gS^{-1}(a\3)a\2\ot h\equiv
\sum v\ot gS^{-1}(a\2)\ot a\1h
$$
modulo $R$, and $R\sbs R'$ since
$$
v\ot g\ot ah=\sum v\ot ga\3S^{-1}(a\2)\ot a\1h\equiv\sum va\1\ot ga\2\ot h
$$
modulo $R'$. Hence $R=R'$, and therefore
$$
M\ot_AH\cong(V\ot H\ot H)/R'\cong V\ot_AX
$$
where $X=H\ot H$ regarded as a left $A$-module with respect to the action of 
$A$ defined by the rule
$$
a\triangleright(g\ot h)=\sum gS^{-1}(a\2)\ot a\1h\,,\qquad a\in A,\ \,g,h\in H.
$$
The linear transformation $\xi$ of $H\ot H$ defined by 
$\xi(g\ot h)=\sum S^{-1}(g\2)\ot g\1h$ has the inverse transformation 
$g\ot h\mapsto\sum S(g\1)\ot g\2h$. Since
$$
\xi(ag\ot h)=\sum S^{-1}(g\2)S^{-1}(a\2)\ot a\1g\1h
=a\triangleright\bigl(\xi(g\ot h)\bigr),
$$
$\xi$ gives an isomorphism of $A$-modules $Y\cong X$ where $Y=H\ot H$ with 
the action of $A$ by left multiplications on the first tensorand. It follows 
that
$$
M\ot_AH\cong V\ot_AX\cong V\ot_AY\cong(V\ot_AH)\ot H.
$$
Since $\xi$ is right $H$-linear with respect to the action of $H$ by right 
multiplications on the second tensorand of $H\ot H$, we get an isomorphism of 
right $H$-modules, as stated in Claim 4.

\medbreak
We are ready now to verify flatness of $H$ over $A$. Let $0\to V'\to V\to V''\to0$ 
be an exact sequence of right $A$-modules. It gives rise to an exact sequence of 
right $A$-modules $0\to M'\to M\to M''\to0$ where
$$
M'=V'\ot H,\qquad
M=V\ot H,\qquad
M''=V''\ot H
$$
with the action of $A$ as specified earlier in the case of $M$. Since the right 
quotient ring $Q(A)$ is left $A$-flat, the sequence
$$
0\to M'\ot_AQ(A)\to M\ot_AQ(A)\to M''\ot_AQ(A)\to0
$$
is exact as well. This sequence of right $Q(A)$-modules splits since all 
terms in it are projective by Claim 3. Applying the functor 
$?\ot_{Q(A)}Q(H)$, we get an exact sequence
$$
0\to M'\ot_AQ(H)\to M\ot_AQ(H)\to M''\ot_AQ(H)\to0.
$$
Note that $M\ot_AQ(H)\cong(M\ot_AH)\ot_HQ(H)\cong(V\ot_AH)\ot Q(H)$ in view of 
Claim 4, and this isomorphism is functorial in $V$. Therefore the previous 
exact sequence can be rewritten as
$$
0\to(V'\ot_AH)\ot Q(H)\to(V\ot_AH)\ot Q(H)\to(V''\ot_AH)\ot Q(H)\to0.
$$
Since the final tensoring in all terms here is performed over the ground 
field, we deduce that the sequence $0\to V'\ot_AH\to V\ot_AH\to V''\ot_AH\to0$ 
is exact. Thus the functor $?\ot_AH$ on the category of right $A$-modules is 
exact, which means that $H$ is indeed left $A$-flat.

There remains the question of faithful flatness. Dealing with it is based on 
the following observation:

\proclaim
Claim 5.
If $V\ne0,$ but $V\ot_AH=0,$ then $A$ has an $\Hd$-stable right ideal $I$ 
such that $I\ne0$ and $I\ne A$.
\endproclaim

Let $M$ be as defined earlier. Applying Claim 4, we obtain $M\ot_AH=0$. Hence 
$M\ot_AQ(H)=0$ too, and so $F\ot_{Q(A)}Q(H)=0$ where we put $F=M\ot_AQ(A)$. 
This is only possible when $F=0$ since $F^n$ is a free $Q(A)$-module for some 
$n>0$, according to Claim 3. Thus $M\ot_AQ(A)=0$. This means that any finite 
subset of $M$ is annihilated by a regular element of $A$. 

Now recall that $M$ is an $\Hd$-equivariant right $A$-module and the action of 
$\Hd$ on $M$ is locally finite. Since $M\ne0$, there exists a finite dimensional 
$\Hd$-submodule $0\ne U\sbs M$. Put $I=\{a\in A\mid Ua=0\}$. Then $I$ is an 
$\Hd$-stable right ideal of $A$. It contains a regular element of $A$ since $I$ 
coincides with the annihilator of any finite basis of $U$. Hence $I\ne0$. On 
the other hand, $1\notin I$ since $U\ne0$. Thus Claim 5 has been verified.

\medskip
Suppose now that $A$ is a Hopf subalgebra. Recall that the $\Hd$-submodules of 
$H$ are precisely the right coideals of $H$. Since $\De(A)\sbs A\ot A$, the 
$\Hd$-submodules of $A$ are precisely the right coideals of $A$. By the Hopf 
module argument recalled in the proof of Lemma 4.2 any Hopf algebra contains 
no nontrivial right ideals which are simultaneously right coideals. This means 
that $A$ has no $\Hd$-stable right ideals other than 0 and the whole $A$. Then, 
by Claim 5, $V\ot_AH\ne0$ for each nonzero right $A$-module $V$. We have seen 
already that $H$ is left $A$-flat, and so we obtain faithful flatness on the 
left.

We have mentioned that $S:H\to H$ is bijective. Since $A$ satisfies the same 
assumptions as $H$, its antipode is also bijective. But the antipode of $A$ is 
the restriction of $S$ to $A$. Thus $S$ is an antiautomorphism of $H$ mapping 
$A$ onto itself. From this it is clear that $H$ is faithfully $A$-flat on both 
sides.
\endproof

\proclaim
Theorem 4.5.
Let $A$ be a right Noetherian right coideal subalgebra of a residually finite 
dimensional Noetherian Hopf algebra $H$. Then $A$ has a right Artinian classical 
right quotient ring{\rm,} and $H$ is left $A$-flat. Moreover{\rm,} if $A$ is a 
Hopf subalgebra{\rm,} then $H$ is left and right faithfully $A$-flat.
\endproclaim

\Proof.
By Lemma 4.2 $A$ is an $\Hd$-prime $\Hd$-module algebra. Since the action of 
$\Hd$ on $A$ is locally finite, we can apply Theorem 4.1 and conclude that 
$A$ has a right Artinian classical right quotient ring $Q(A)$. The same result 
can be applied to $H$ in place of $A$.

Thus we meet the hypothesis of Theorem 4.4, and everything follows.
\endproof

\proclaim
Corollary 4.6.
Retain all assumptions of Theorem\/ {\rm4.5}. If $A$ is a Hopf subalgebra{\rm,} 
then $H$ is a projective generator in the categories of right and left $A$-modules.
\endproclaim

\Proof.
This follows from a result of Schneider \cite{Schn92, Cor. 1.8}. See also 
Masuoka and Wigner \cite{Ma-W94, Th. 2.1}.
\endproof

\references
\nextref An92
\auth{A.Z.,Anan'in}
\paper{Representability of noetherian finitely generated algebras}
\journal{Arch. Math.}
\Vol{59}
\Year{1992}
\Pages{1-5}

\nextref Coh86
\auth{M.,Cohen}
\paper{Smash products, inner actions and quotient rings}
\journal{Pacific J. Math.}
\Vol{125}
\Year{1986}
\Pages{45-66}

\nextref Goo
\auth{K.R.,Goodearl}
\book{Nonsingular Rings and Modules}
\publisher{Marcel Dekker}
\Year{1976}

\nextref Goo-W
\auth{K.R.,Goodearl;R.B.,Warfield Jr.}
\book{An Introduction to Noncommutative Noetherian Rings}
Second edition,
\publisher{Cambridge Univ. Press}
\Year{2004}

\nextref Ma91
\auth{A.,Masuoka}
\paper{On Hopf algebras with cocommutative coradicals}
\journal{J.~Algebra}
\Vol{144}
\Year{1991}
\Pages{451-466}

\nextref Ma92
\auth{A.,Masuoka}
\paper{Freeness of Hopf algebras over coideal subalgebras}
\journal{Comm. Algebra}
\Vol{20}
\Year{1992}
\Pages{1353-1373}

\nextref Ma-W94
\auth{A.,Masuoka;D.,Wigner}
\paper{Faithful flatness of Hopf algebras}
\journal{J.~Algebra}
\Vol{170}
\Year{1994}
\Pages{156-164}

\nextref Mc-R
\auth{J.C.,McConnell;J.C.,Robson}
\book{Noncommutative Noetherian Rings}
Revised edition,
\publisher{Amer. Math. Soc.}
\Year{2001}

\nextref Mo
\auth{S.,Montgomery}
\book{Hopf Algebras and their Actions on Rings}
\publisher{Amer. Math. Soc.}
\Year{1993}

\nextref Mo-Sch95
\auth{S.,Montgomery;H.-J.,Schneider}
\paper{Hopf crossed products, rings of quotients, and prime ideals}
\journal{Adv. Math.}
\Vol{112}
\Year{1995}
\Pages{1-55}

\nextref Nich-Z89
\auth{W.D.,Nichols;M.B.,Zoeller}
\paper{A Hopf algebra freeness theorem}
\journal{Amer. J. Math.}
\Vol{111}
\Year{1989}
\Pages{381-385}

\nextref Rad77
\auth{D.E.,Radford}
\paper{Pointed Hopf algebras are free over Hopf subalgebras}
\journal{J.~Algebra}
\Vol{45}
\Year{1977}
\Pages{266-273}

\nextref Scha00
\auth{P.,Schauenburg}
\paper{Faithful flatness over Hopf subalgebras: counterexamples}
\InBook{Interactions between ring theory and representations of algebras}
\publisher{Marcel Dekker}
\Year{2000}
\Pages{331-344}

\nextref Schn92
\auth{H.-J.,Schneider}
\paper{Normal basis and transitivity of crossed products for Hopf algebras}
\journal{J.~Algebra}
\Vol{152}
\Year{1992}
\Pages{289-312}

\nextref Schn93
\auth{H.-J.,Schneider}
\paper{Some remarks on exact sequences of quantum groups}
\journal{Comm. Algebra}
\Vol{21}
\Year{1993}
\Pages{3337-3357}

\nextref Sk06
\auth{S.,Skryabin}
\paper{New results on the bijectivity of antipode of a Hopf algebra}
\journal{J.~Algebra}
\Vol{306}
\Year{2006}
\Pages{622-633}

\nextref Sk07
\auth{S.,Skryabin}
\paper{Projectivity and freeness over comodule algebras}
\journal{Trans. Amer. Math. Soc.}
\Vol{359}
\Year{2007}
\Pages{2597-2623}

\nextref Sk08
\auth{S.,Skryabin}
\paper{Projectivity of Hopf algebras over subalgebras with semilocal central localizations}
\journal{J.~$K$-Theory}
\Vol{2}
\Year{2008}
\Pages{1-40}

\nextref Sk10
\auth{S.,Skryabin}
\paper{Models of quasiprojective homogeneous spaces for Hopf algebras}
\journal{J.~Reine Angew. Math.}
\Vol{643}
\Year{2010}
\Pages{201-236}

\nextref Sk11
\auth{S.,Skryabin}
\paper{Structure of $H$-semiprime Artinian algebras}
\journal{Algebr. Represent. Theory}
\Vol{14}
\Year{2011}
\Pages{803-822}

\nextref Sk-Oy06
\auth{S.,Skryabin;F.,Van Oystaeyen}
\paper{The Goldie theorem for $H$-semiprime algebras}
\journal{J.~Algebra}
\Vol{305}
\Year{2006}
\Pages{292-320}

\nextref St
\auth{B.,Stenstr\"om}
\book{Rings of Quotients}
\publisher{Springer}
\Year{1975}

\nextref Tak79
\auth{M.,Takeuchi}
\paper{Relative Hopf modules---equivalences and freeness criteria}
\journal{J.~Algebra}
\Vol{60}
\Year{1979}
\Pages{452-471}

\nextref Wu-Zh03
\auth{Q.-S.,Wu;J.J.,Zhang}
\paper{Noetherian PI Hopf algebras are Gorenstein}
\journal{Trans. Amer. Math. Soc.}
\Vol{355}
\Year{2003}
\Pages{1043-1066}

\endreferences
\bye